\documentclass[12pt,a4paper]{amsart}

\usepackage{amsmath}
\usepackage{amssymb}
\usepackage{latexsym}
\usepackage{xypic}
\usepackage{amscd,amsthm}

\usepackage{graphicx}
\usepackage{psfrag}
\usepackage{epic}
\usepackage{eepic}

\def\P{{\mathbb P}}

\def\R{{\mathbb R}}
\def\Z{{\mathbb Z}}
\def\C{{\mathbb C}}

\newtheorem*{theorem*}{Theorem}
\newtheorem*{corollary*}{Corollary}
\newtheorem{theorem}{Theorem}[section]
\newtheorem{definition}[theorem]{Definition}

\newtheorem{example}[theorem]{Example}
\newtheorem{corollary}[theorem]{Corollary}
\newtheorem{proposition}[theorem]{Proposition}

\begin{document}
\title[Conifold degenerations of Fano $3$-folds]{Conifold
degenerations of Fano $3$-folds as hypersurfaces
in toric varieties}

\author[Victor Batyrev]{Victor Batyrev}
\address{Mathematisches Institut, Universit\"at T\"ubingen,
Auf der Morgenstelle 10, 72076 T\"ubingen, Germany}
\email{ victor.batyrev@uni-tuebingen.de}

\author[Maximilian Kreuzer]{Maximilian Kreuzer}
\address{Institut f\"ur Theoretische Physik,
Technische Universit\"at Wien,  Wiedner Hauptstr. 8-10/136
A-1040 Vienna, Austria}
\email{Maximilian.Kreuzer@tuwien.ac.at}


\begin{abstract}
There exist exactly $166$
$4$-dimensional reflexive polytopes $\Delta$
such that the corresponding $4$-dimensional
Gorenstein toric Fano variety $\P_\Delta$ has
at worst terminal singularities in codimension
$3$ and the anticanonical divisor of  $\P_\Delta$ is divisible
by $2$ in its Picard group.
For every such a polytope $\Delta$, one naturally obtains
a family  ${\mathcal F}(\Delta)$ of Fano
hypersurfaces  $X \subset \P_\Delta$
with
at worst conifold singularities. A generic
$3$-dimensional Fano hypersurface $X \in {\mathcal F}(\Delta)$
can be interpreted  as a flat conifold degeneration of
some smooth Fano $3$-folds $Y$ whose classification up to deformation
was obtained by Iskovskikh, Mori and Mukai.
In this case, both Fano varieties $X$ and $Y$ have the same
Picard number $r$. Using toric mirror symmetry, we define
a $r$-dimensional generalized hypergeometric power series
$\Phi$  associated to the dual reflexive polytope
$\Delta^*$.
We show that if $r =1$ then $\Phi$
is a normalized regular solution of a modular
$D3$-equation that appears in the Golyshev correspondence.
We expect that the power series $\Phi$
can be used to compute the small quantum cohomology
ring of all Fano $3$-folds $Y$ with the Picard number $r \geq 2$
if $Y$ admit a conifold degeneration $X \in {\mathcal F}(\Delta)$.
\vspace{-0.5ex}\end{abstract}

\maketitle

\section*{Introduction}

A smooth projective $d$-dimensional algebraic variety
$V$ is called {\em Fano $d$-fold} if the anticanonical class
$-K_V$ is an ample Cartier divisor. The maximal integer
$m$ such that $-K_V = mL$ for some Cartier divisor on $V$
is called the {\em index of $V$}.

It has been proved by Kollar, Mori and Myaoka that there exist
only finitely many smooth Fano $d$-folds of fixed dimension $d$ up
to deformation \cite{KMM92}. The classification of smooth Fano $d$-folds
is known only for  $d \leq 3$. There exist $17$ types of Fano $3$-folds with
the Picard number $1$ (they were classified by Iskovskikh
\cite{Isk77,Isk78}) and $105$ types of Fano $3$-folds with
the Picard number
$\geq 2$ (see the classification due to Mori and Mukai \cite{MoMu81,MoMu03}).
Among all $122$ types of Fano $3$-folds there are  exactly
$18$ Fano $3$-folds which are toric varieties (these varieties
 were classified
by the first author \cite{Ba81} and Watanabe-Watanabe
\cite{WaWa82}).

It is known that for any fixed dimension $d$ there exist only finitely
many $d$-dimensional toric Fano varieties with at worst
Gorenstein singularities \cite{Ba82}. These varieties
are of  particular  interest because of their relation to
mirror symmetry \cite{Ba94}. We also remark that there is a bijection
between
$d$-dimensional Gorenstein toric Fano varieties
ans so called {\em reflexive} $d$-dimensional lattice
polytopes.
The second author   together with Skarke has classified
all reflexive polytopes of dimension $3$ and $4$ \cite{KS98,KS00}.
This classification has been obtained by a computer and it
contains of $4319$ $3$-dimensional reflexive polytopes
and $473\, 800 \, 776$  $4$-dimensional
reflexive polytopes.

In this paper we  connect the
Iskovskikh-Mori-Mukai
classification of Fano $3$-folds to  the classification of
Fano hypersurfaces with conifold singularities
in  $4$-dimensional Gorenstein toric Fano varieties $\P_\Delta$
of index $2$
corresponding to some special $4$-dimensional reflexive
polytopes $\Delta$. Our approach is motivated by
a theorem of  Namikawa \cite{Na97} that claims that
every 3-dimensional
Gorenstein Fano variety $X$ with at worst terminal
singularities admits a flat smoothing to a Fano $3$-fold $Y$.
In this paper
we construct many examples of such singular  3-dimensional
Gorenstein Fano varieties as generic hypersurfaces in
$4$-dimensional Gorenstein toric Fano varieties $\P_\Delta$ corresponding
to  $4$-dimensional
reflexive polytopes $\Delta$ that satisfy the following two
conditions:

\begin{enumerate}
\item $\Delta$ is isomorphic to $2\Delta'$ for some
$4$-dimensional lattice polytope $\Delta'$ (i.e.,
$\Delta$ is "divisible by $2$");

\item all $2$-dimensional faces of the dual
reflexive polytope $\Delta^*$ are either basic
triangles, or parallelograms consisting
of two basic triangles.
\end{enumerate}

Among $473\, 800 \, 776$  $4$-dimensional
reflexive polytopes  we find  $5363$
polytopes satisfying the
 first condition and $198\, 849$ polytopes
satisfying the second one.
By intersecting these two sets, we obtain
$166$ $4$-dimensional reflexive polytopes $\Delta$ satisfying
both conditions $(1)$ and $(2)$.
Among $166$ polytopes there are exactly $100$ reflexive polytopes
$\Delta \cong 2 \Delta'$ such that $\Delta'$
that are pyramids over $3$-dimensional reflexive
polytopes $\Theta$ that define
$3$-dimensional toric Fano varieties $\P_\Theta$ with
at worst conifold singularities (i.e. $\P_\Theta$ is
a conifold toric degeneration of a  smooth Fano $3$-fold).
We remark that these toric degenerations were classified  independently
by Nill \cite{Ni05} and Galkin \cite{Ga08,Ga12}.
The interest to toric degenerations is motivated
by the mirror symmetry for complete
intersections in Fano varieties \cite{BCKS98,BCKS00,Ba04}.
On the other hand, toric degenerations are useful
for constructing mirros of Fano varieties via
Landau-Ginzburg models \cite{Pr07c}.

\section{Fano hypersurfaces}

Let $M \cong \Z^d$ be a free abelian group of rank $d$ and
$N := {\rm Hom}(M, \Z)$  be the dual group with the natural
pairing $\langle *, * \rangle \; : \; M \times N  \to \Z$.
We consider also $\R$-vector spaces $M_\R := M \otimes_\Z \R$
and  $N_\R := N \otimes_\Z \R$ together with the pairing
 $\langle *, * \rangle \; : \; M_\R \times N_\R  \to \R$.
Recall that a $d$-dimensional polytope
 $\Delta \subset M_\R$ is called {\em reflexive} if
it contains $0 \in M$ in its interior, all vertices of $\Delta$
belong to the lattice $M \subset M_\R$,
and all vertices of the dual polytope
\[ \Delta^* := \{ y \in N_\R\; | \; \langle x, y \rangle \geq -1
\;\; \forall x \in \Delta \} \]
belong to the dual lattice $N \subset N_\R$. If $\Delta$ is reflexive then
$\Delta^*$ is also reflexive and $(\Delta^*)^* = \Delta$.
If $\Theta \subset \Delta$ is a proper face of a reflexive pooytope
$\Delta$ then we set
$$\Theta^*:= \{  y \in \Delta^*\; | \; \langle x, y \rangle = -1
\;\; \forall x \in \Theta \}$$
and call $\Theta^*$ the dual face of $\Delta$. One has
$\dim \Theta + \dim \Theta^* = d-1$.
If $\Delta \subset M_\R$ is a $d$-dimensional reflexive polytope then
we denote by $\Sigma(\Delta^*)$ (or simply by $\Sigma$)
the complete fan consisting of all cones
$\sigma(\Theta^*):= \R_{\geq 0}\Theta^*  \subset N_\R$ where $\Theta^*$
runs  over all proper faces $\Theta^*$
of the dual reflexive polytope $\Delta^*$. The fan $\Sigma(\Delta^*)$
defines a $d$-dimensional projective variety over $\C$ that we denote by
$\P_\Delta$. One can define the projective toric variety
$\P_\Delta$ also as ${\rm Proj}$ of  the graded semigroup algebra over
$\C$ corresponding to the graded monoid of all lattice points $(k,
{\bf m}) \in \Z_{\geq 0}  \times  M$ such that ${\bf m} \in k \Delta$.
There exists a natural stratification of $\P_\Delta$ by torus orbits
$T_\Theta$ where $\Theta$ runs over all faces of $\Delta$. Then
$T = T_\Delta$ is the open dense orbit and $\dim T_\Theta = \dim \Theta$
$\forall \Theta \subset \Delta$. A generic ample divisor $X \subset \P_\Delta$
admit a stratification
\[ X = \bigcup_{\dim \Theta >0} X_\Theta \]
where $X_\Theta := X \cap T_\Theta$. The toroidal singularities of
$X$ along $X_\Theta$ are determined by the cone $\sigma(\Theta^*)$.
By the combinatorial characterisation of terminal singularities,
one obtains that the Gorenstein singularities of $X$ along $X_\Theta$ are
terminal if and only if all lattice points in the
dual face $\Theta^* \subset \Delta^*$ are its vertices.
It is no difficult to show that there exist exactly two
$2$-dimensional polygons $P \subset \R^2$ with vertices in $\Z^2$
(up to an $\Z$-isomorphism) such that $P \cap \Z^2$ is exactly the set
of vertices of $P$:

\begin{itemize}

\item  the triangle with vertices $(0,0)$, $(1, 0)$, $(0,1)$;

\item  the unit square with vertices  $(0,0)$, $(1, 0)$, $(0,1)$,
$(1,1)$.

\end{itemize}

This implies that Gorenstein singularities in codimension $3$
are locally conifold singularities defined by the equation:
\[ xy = zt   \sim_\C x^2 + y^2 + z^2 + t^2 = 0. \]

We remark that the anticanonical class $-K_{\P_\Delta}$ is divisible
in ${\rm Pic}(\P_\Delta)$ by some positive integer $k$ (i.e.$-K_{\P_\Delta} =
k H$ for some $H \in {\rm Pic}(\P_\Delta)$)
if and only if the reflexive polytope $\Delta$ is isomorphic
to $k\Delta'$ for some lattice polytope $\Delta' \subset M_\R$.

From now on we consider the case $d =4$.

\begin{definition}
Let $\Delta$ be a $4$-dimensional
reflexive polytope such that $\Delta \cong 2 \Delta'$ for some lattice
polytope $\Delta' \subset M_\R$. We denote by ${\mathcal F}(\Delta)$
the family of Fano hypersufaces in $\P_\Delta$ such that the Newton polytope
of their equations equals $\Delta'$.
\end{definition}

\begin{proposition}
A generic Fano hypersuface $X \in {\mathcal F}(\Delta)$ hast at worst
Gorenstein terminal singularities if and only if every $2$-dimensional
face $\Theta^* \subset \Delta^*$ is a lattice polygon that is isomorphic
to the standard triangle, or to the unit square.
\label{sing}
\end{proposition}

Using the classification of all $4$-dimensional reflexive polytopes
and \ref{sing}, we obtain

\begin{theorem}
There exist exactly $166$ $4$-dimensional reflexive polytopes
$\Delta \cong 2 \Delta'$ such that
a generic Fano hypersurface $X \in {\mathcal F}(\Delta)$
in  $4$-dimensional Gorenstein toric Fano variety
$\P_\Delta$ has at worst terminal Gorenstein singularities
\end{theorem}

On the other hand, there is a following statement
due to Namikawa \cite[Theorem 11]{Na97}:

\begin{theorem}
Let $X$ be $3$-dimensional Fano variety
with Gorenstein terminal singularities.
Then there is a flat deformation of $X$ to a smooth
Fano $3$-fold $Y$.
\end{theorem}

\begin{corollary}
Every $4$-dimensional reflexive polytope $\Delta$ satisfying
the conditions $(1)$ and $(2)$
gives rise to a deformation type of smooth Fano $3$-folds
$Y$ in the Iskovskikh-Mori-Mukai classification.
\end{corollary}

Our next pupose will be to explain how to compute the
the topological invariants of the smoothing $Y$
via combinatorics of reflexive polytopes $\Delta$ and
$\Delta^*$.

Since the degree of projective varieties
remains unchanged under flat deformations
we get:

\begin{proposition}
The anticanonical degree $(-K_Y)^3$ of the smoothing
$Y$ of Fano hypersurfaces $X \in {\mathcal F}(\Delta)$ is
equal to
\[  4! {\rm vol}(\Delta') = \frac{4!{\rm vol}(\Delta)}{16}. \]
\end{proposition}

Our next interest is the Picard number of $Y$. It is
important observation that in our situation also  the Picard group 
 remains unchanged
under the flat deformation \cite{JR06}.  Therefore, it remains
to compute the Picard number of the Fano hypersurface $X \subset
\P_\Delta$. By Lefschetz-type arguments, the latter ist
equal to the Picard number of the toric variety
$\P_\Delta$.

Let $\Delta \subset M_\R$ be a $4$-dimensional
reflexive polytope satisfying $(1)$ and $(2)$
and let $\{ v_1, \ldots, v_n \}$ be the set of vertices
of the dual reflexive polytope $\Delta^*$.
Denote by $PL(\Delta*)$ the sublattice in $\Z^n$
consisiting of all $n$-tuples $(l_1, \ldots, l_n) \in \Z^n$ such that
$l_{i_1} + l_{i_3} = l_{i_2} + l_{i_4}$  holds
whenever $v_{i_1}, v_{i_2}, v_{i_3}, v_{i_3}$ are vertices of
a $2$-dimensional face $\Theta^*$ of $\Delta^*$ satisfying the
equation $v_{i_1} + v_{i_3} = v_{i_2} + v_{i_4}$.
We set
\[ rk(\Delta^*) := n - {\rm rank}\,  PL(\Delta^*). \]

\begin{proposition}
Consider the monomorphism
\[ \varphi\, : \, M \to \Z^n \]
\[ m \mapsto (\langle m, v_1 \rangle, \ldots, \langle m, v_n \rangle) \]
Then the image of $\varphi$ is contained in $PL(\Delta^*) \subset \Z^n$
and the Picard group of the toric variety $\P_\Delta$ is
isomorphic to the group
\[ PL(\Delta)/\varphi(M). \]
In particular, the Picard number of the toric variety $\P_\Delta$
is equal to
\[ n -4 - rk(\Delta^*), \]
where $n = {\rm vert}(\Delta^*)$ is the number of vertices
of $\Delta^*$.
\end{proposition}

{\em Proof.} The statement follows from the standard computation of the 
Picard group of a toric variety \cite{Ei92}. \hfill $\Box$ 

\begin{corollary}
The second Betti number of a Fano  $3$-fold $Y$ that admits a flat
conifold degeneration to a Fano hypersurface $X \in {\mathcal F}(\Delta)$
is at most $5$. Its distribution for $166$ reflexive
polytopes $\Delta$ is given by the following table:
\[ \begin{tabular}{|c|c|c|c|c|c|}
\hline
$B_2(Y)$ &  1 & 2 & 3 & 4 & 5 \\
\hline
the number of  polytopes $\Delta$ & 23 &  69 & 54 & 18 & 2 \\
\hline
\end{tabular}
\]

\end{corollary}

The last topological invariant that we want to compute is
the Hodge number $h^{2,1}(Y) = 1/2 B_3(Y)$. This number
depends on the number of conifold singularities on the
Fano hypersurface $X \in {\mathcal F}(\Delta)$.

\begin{proposition}
Let
$$sq(\Delta^*):= \sum_{\Theta^* \subset \Delta, \; \dim \Theta^* =2}
(| \Theta^* \cap N| -3)  $$
be the number of $2$-dimensional faces $\Theta^* \subset \Delta^*$
containing $4$ vertices (i.e. the number of ``squares'').
Then the number
of conifold points in $X$ is equal to
\[  dp(\Delta^*):= \sum_{\Theta^* \subset \Delta, \; \dim \Theta^* =2}
(| \Theta^* \cap N| -3)(|\Theta \cap M| -1). \]
\end{proposition}

There exists
a small resolution $\widehat{X} \to X$
of these singularities such that
$\rho(\widehat{X}) - \rho({X}) =
rk(\Delta^*)$ where $\rho(V)$ denotes the Picard number of $V$.

\begin{definition}
Let us put $py(\Delta) :=1$ if $\Delta = 2\Delta'$
and $\Delta'$ is a pyramid over $3$-dimensional reflexive polytope and
$py(\Delta) :=0$ otherwise.
\end{definition}

\begin{proposition}
\[ h^{2,1}(Y)  = \frac{1}{2} B_3(Y) = 1 + dp(\Delta^*)   -rk(\Delta^*) -
py(\Delta) = 1 - py(\Delta) + dp(\Delta^*) + \rho({X}) -
\rho(\widehat{X}). \]
\end{proposition}

{\em Proof.} The smooth Fano $3$-fold $Y$ is obtained from the smooth almost 
Fano $3$-fold $\widehat{X}$ by so called {\em conifold transition}. 
First we compute the Hodge number $h^{2,1}(X)$ of the Fano hypersurface 
$X \subset \P_\Delta$ using the formulas of Danilov-Khovanskii \cite{DKh86}. 
We note that the Newton polytope of the equation of $X$ is $\Delta'$. 
There is no interior points in $\Delta'$ and there is exactly one 
interion point in $2\Delta' = \Delta$. A codimension $1$ face of $\Delta'$ contains
an interior lattice point if and only if $\Delta'$ is a 
pyramid over a $3$-dimensional 
reflexive polytope. In this case only one facet has an 
interior lattice point. 
This shows that $h^{2,1}(X) = 1 - py(\Delta)$. Now we 
apply standard Clemens arguments
to get $$h^{2,1}(Y) = h^{2,1}(X) +  dp(\Delta^*)  + \rho({X}) -
\rho(\widehat{X}).$$
\hfill $\Box$

\section{The power series $\Phi$}

Let $\Delta$ be a $4$-dimensional reflexiv polytopes satisfying the conditions
(1) and (2) as above. Denote by $n := {\rm vert}(\Delta^*)$ the number
of vertices of the dual polytope $\Delta^*$.

Consider the
lattice of rank $n- 4$:
\[  \Lambda(\Delta^*):= \{ {\bf k} =
(k_1, \ldots, k_n) \in \Z^n  \; |\;
\sum_{i =1}^n k_i v_i = 0 \}. \]
One has the natural pairing
\[  PL(\Delta^*)\times \Lambda(\Delta^*)  \to \Z \]
\[ ({\bf l}, {\bf k}) \mapsto \sum_{i =1}^n l_i k_i \]
that vanish for all ${\bf l} = (l_1, \ldots, l_n) \in \varphi(M)$, because
for any $m \in M$ and any   ${\bf k} = (k_1, \ldots, k_n) \in
\Lambda(\Delta^*)$ one has
\[ \sum_{i=1}^n k_i \langle m, v_i \rangle  =  \langle m, \sum_{i=1}^n k_i v_i
\rangle =  \langle m, 0 \rangle =0. \]
So one obtains  the pairing
between ${\rm Pic}(\P_\Delta) \cong  PL(\Delta^*)/\varphi(M)$
and  $\Lambda(\Delta^*)$. Let $\lambda_1, \ldots, \lambda_r \in
 PL(\Delta^*)/\varphi(M)$ be a $\Z$-basis.  We denote by  $\Lambda^+(\Delta^*)$
the semigroup $ \Lambda(\Delta^*) \cap \Z_{\geq 0}^n$ and by
$\kappa$ the element of  $PL(\Delta^*)/\varphi(M) \cong {\rm Pic}(\P_\Delta)
$ such that
$2\kappa$ ist the anticanonical class of $\P_\Delta$, i. e., 
$k_1 + \cdots +k_n = 2(\kappa, {\bf k})$. 
Define the multidimensional series $\Phi$ by the formula
\[ \Phi(t_1, \ldots, t_r)  =
\sum_{ {\bf k} \in  \Lambda^+(\Delta^*) } \frac{((\kappa, {\bf k})!)^2}
{k_1! \cdots k_n!} {t_1}^{(\lambda_1, {\bf k})} \cdots
{t_r}^{(\lambda_r, {\bf k})}. \]
These series were first
suggested in \cite{BaSt95} in connection to the mirror symmetry
for complete intersections in toric varieties. In our situation, 
it corresponds to mirrors of $K3$-surfaces that are complete intersections 
of two divisors in the $4$-dimensional Gorenstein toric Fano variety 
$\P_\Delta$. 

There is a specialization $\Phi_0(t)$ of $\Phi$ to a 
$1$-parameter series corresponding 
to the class $\kappa \in {\rm Pic}(\P_\Delta)$ that 
restricts to the anticanonical 
class in ${\rm Pic}(X)$:
\[ \Phi_0(t) := \sum_{ {\bf k} \in  \Lambda^+(\Delta^*) } 
\frac{((\kappa, {\bf k})!)^2}
{k_1! \cdots k_n!} {t}^{(\kappa, {\bf k})}. \]

\begin{example} 
{\rm There exists only one simplex $\Delta$ 
among 166 reflexive polytopes satisfying the conditions 
$(1)$ and $(2)$.  The vertices of the dual reflexive simplex 
$\Delta^*$ satisfy the 
single relation 
\[ 4v_0 + v_1 + v_2 + v_3 + v_4 =0. \]
The corresponding hypersurface $X$ is isomorphic to $\P^3$ that considered as 
a hypersurface of degree $4$ in the $4$-dimensional 
weighted projective space $\P(4,1,1,1,1)$. So we have 
\[ \Phi_0(t) = \sum_{k \geq 0} \frac{(({4k})!)^2}
{(4k)! (k!)^4} {t}^{4k} = \sum_{k \geq 0} \frac{({4k})!}
{(k!)^4} {t}^{4k}. \]
} 
\end{example} 

\begin{example} 
{\rm There exist exactly $2$ reflexive $4$-dimensional polytopes from the list 
of $166$ ones such that $\Delta^*$ has $6$ vertices satisfying two independent 
relations 
\[ 3v_0 + v_1 + v_3 + v_4 =0, \;\; 3v_0 + v_1 + v_2 + v_5 =0 \]
with the corresponding power series 
\[  \Phi_0(t) = \sum_{k \geq 0} \sum_{k_1 + k_2 = k} 
\frac{(({3k})!)^2}{(3k)! (k!)(k_1!)^2(k_2!)^2} t^{3k} = \sum_{k \geq 0} 
\frac{({3k})!(2k)!}{ (k!)^5} t^{3k}, \] 
or the relations 
\[ v_0 + v_1 + v_3 + v_4 =0, \;\; v_0 + v_1 + v_2 + v_5 =0 \]
with the corresponding power series
\[  \Phi_0(t) = \sum_{k \geq 0} \sum_{k_1 + k_2 = k}
\frac{(({2k})!)^2}{((k)!)^2(k_1!)^2(k_2!)^2} t^{2k} = \sum_{k \geq 0}
\frac{((2k)!)^3}{ (k!)^6} t^{2k}. \]
In the above calculations we used the equaity 
\[ \sum_{k_1 + k_2 = k} \frac{(k!)^2}{(k_1!)^2 (k_2!)^2} = 
\frac{(2k)!}{(k!)^2}. \]}
\end{example}

\section{Hypersurfaces with the Picard number  $1$}

There exists exactly $23$ $4$-dimensional reflexive 
polytopes $\Delta$ satisfying the conditions $(1)$ and $(2)$ 
such that a generic Fano hypersurface $X \in {\mathcal F}(\Delta)$ 
has the Picard number $1$. 

We will use the following notations:  $\deg :=( -K_Y)^3$, $h^{2,1} := 
h^{2,1}(Y)$, $rk:= rk(\Delta^*)$, $dp:= dp(\Delta^*)$, 
$sq:= sq(\Delta^*)$, $py:= py(\Delta)$, ${\rm vert}(\Delta^*)$ 
denotes the number of vertices of $\Delta^*$.  The deformation 
type of a Fano $3$-fold $Y$ with the Picard number 
$1$ is completely determined by the index $m$ and 
the integer $\deg/(2m^2)$. The pair $(m, \deg/(2m^2))$ is called
the type of the Fano $3$-fold $Y$.  There exist exactly $17$ types 
of smooth Fano $3$-folds with the Picard number $1$: 
\[ (4,2), (3,3) \]
\[ (2,k), \;\; \; 1 \leq k \leq 5 \] 
\[ (1,k), \; \; \; 1 \leq k \leq 9, \;\; k=11\]
Among $17$ types there are $13$ ones that admit  conifold degenerations  
$X \in {\mathcal F}(\Delta)$. There remaining $4$ types 
are hypersurfaces or complete intersections in weighted 
projective spaces:
\[  (2,1) \; :  \; V_6 \subset \P(3,2,1,1,1) \]
\[ (1,1)\;  :\;  V_6 \subset \P(3,1,1,1,1) \]
\[ (1, 2) \; : \; V_4 \subset \P^4 \]
\[ (1,3) \; : \; V_{2,3} \subset \P^5 \]
 
The table below describes vertices of $23$ dual polytopes 
$\Delta^*$ together with its properties and 
topological  invariants of $Y$. 

\newpage

{\tiny

$----------------------------------$

Type $(4,2)$

$----------------------------------$

$V(1)$: 
$\deg = 64$, $h^{2,1} =0$, $rk =0$, $dp =0$, $sq =0$, $py =1$, 
${\rm vert}(\Delta^*) =5$

\begin{tabular}{|c|c|c|c|c|}
1& 0& 0& 0& -4  \\
0& 1& 0& 0& -1 \\
0& 0& 1& 0& -1  \\
0& 0& 0& 1& -1
\end{tabular}

$----------------------------------$

Type $(3,3)$

$----------------------------------$

$V(2)$: 
$\deg = 54$, $h^{2,1} =0$, $rk =1$, $dp =1$, $py =1$, 
${\rm vert}(\Delta^*) =6$:

\begin{tabular}{|c|c|c|c|c|c|}
1& 0& 0& 0& -3 & -3 \\
0& 1& 0& 0& -1 & -1\\
0& 0& 1& 0& 0  & -1 \\
0& 0& 0& 1& -1 & 0
\end{tabular}

$----------------------------------$


Type $(2,2)$

$----------------------------------$

$V(3)$: $\deg = 16$, $h^{2,1} =10$, $rk =3$, $sq =6$, $dp =12$, $py =0$, 
${\rm vert}(\Delta^*) = 8$:

\begin{tabular}{|c|c|c|c|c|c|c|c|}
  1&-1&-1 &-1&-1  & 1 & 1 &-1  \\
  1& 0& 2 & 0& 0  &-1 &-1 & 1 \\
 -1& 0&-1 & 1& 0  & 0 & 1 & 0 \\
 -1& 1&-1 & 0& 0  & 1 & 0 & 0
\end{tabular}

$----------------------------------$

Type $(2,3)$

$----------------------------------$

$V(4)$: $\deg = 24$, $h^{2,1} =5$, $rk =2$, $sq =3$, $dp =6$, $py =0$, 
${\rm vert}(\Delta^*) = 7$:

\begin{tabular}{|c|c|c|c|c|c|c|}
 -1&-1& 1 &-1&-1   & 1 &-1  \\
  0& 2&-1 & 0& 0   &-1 & 1 \\
  0&-1& 0 & 1& 0   & 1 & 0 \\
  1&-1& 1 & 0& 0   & 0 & 0
\end{tabular}

$----------------------------------$

Type $(2,4)$

$----------------------------------$

$V(5)$: 
$\deg = 32$, $h^{2,1} =2$, $rk =1$, $sq =1$, $dp =2$, $py =0$, 
${\rm vert}(\Delta^*) =6$

\begin{tabular}{|c|c|c|c|c|c|}
1& 0& 0& 0& -1 & -1 \\
0& 1& 0& 0& -1 & -1\\
0& 0& 1& 0& 0  & -1 \\
0& 0& 0& 1& -1 & 0
\end{tabular}

$----------------------------------$

$V(6)$: $\deg = 32$, $h^{2,1} =2$, $rk =4$, $sq =6$, $dp =6$, $py =1$, 
${\rm vert}(\Delta^*) = 9$:

\begin{tabular}{|c|c|c|c|c|c|c|c|c|}
-1 &-1& 1&-1 &-1&-1  &-1 &-1 &-1  \\
 0 & 2&-1& 0 & 0& 0  & 2 & 2 & 2 \\
 1 &-1& 0& 0 & 1& 0  & 0 &-1 & 0 \\
 1 &-1& 0& 1 & 0& 0  & 0 & 0 &-1
\end{tabular}

$----------------------------------$

Type  $(2,5)$

$----------------------------------$

$V(7)$: $\deg = 40$, $h^{2,1} =0$, $rk =3$, $sq =3$, $dp =3$, $py =1$, 
${\rm vert}(\Delta^*) = 8$:

\begin{tabular}{|c|c|c|c|c|c|c|c|}
 0&-4&-2 & 1& 0  & 0 &-2 &-2  \\
 0&-1&-1 & 0& 0  & 1 & 0 & 0 \\
 0&-1& 0 & 0& 1  & 0 &-1 & 0 \\
 1&-1& 0 & 0& 0  & 0 & 0 &-1
\end{tabular}

$----------------------------------$

\newpage

Type $(1,4)$:

$----------------------------------$

$V(8)$: $\deg = 8$, $h^{2,1} =14$, $rk =11$, $sq =24$, $dp =24$, $py =0$, 
${\rm vert}(\Delta^*) = 16$:

\begin{tabular}{|c|c|c|c|c|c|c|c|c|c|c|c|c|c|c|c|}
1 & -1 & -1&  1& -1 &  1&  1& -1 & 1 & -1 &-1 & 1 & 1 &-1 & -1 &1  \\
-1&  1 & 0 &  0&  1 & -1&  0&  0 & -1&  1 & 0 & 0 &-1 & 1 &  0 & 0 \\
-1&  1 & 1 &  0&  0 & -1& -1&  0 & 0 &  1 & 1 & 0 & 0 & 0 &  0 &-1 \\
-1&  1 & 1 &  0&  1 &  0&  0&  1 & 0 &  0 & 0 &-1 &-1 & 0 &  0 &-1
\end{tabular}

$----------------------------------$

Type  $(1,5)$:

$----------------------------------$

$V(9)$:  $\deg = 10$, $h^{2,1} =12$, $rk =9$, $sq =18$, $dp =18$, $py =0$, 
${\rm vert}(\Delta^*) = 14$:

\begin{tabular}{|c|c|c|c|c|c|c|c|c|c|c|c|c|c|}
1 & -1 & -1&  1& -1 & -1&  1&  1 & -1& -1 &1 & 1 &-1 &-1 \\
-1&  1 & 0 & -1&  1 &  0&  0& -1 & 1 &  1 &0 &-1 & 0 & 0\\
-1&  1 & 1 &  0&  0 &  0& -1& -1 & 1 &  0 &-1& 0 & 1 & 0\\
-1&  1 & 1 &  0&  1 &  1&  0&  0 & 0 &  0 &-1&-1 & 0 & 0
\end{tabular}

$----------------------------------$

Type $(1,6)$:
 
$----------------------------------$

$V(10)$: $\deg = 12$, $h^{2,1} =7$, $rk =7$, $sq =13$, $dp =13$, $py =0$, 
${\rm vert}(\Delta^*) = 12$:

\begin{tabular}{|c|c|c|c|c|c|c|c|c|c|c|c|}
1  & -1 & -1& -1& 1 &  1& 1  & -1 &-1 &-1 & -1 &1\\
-1 & 0  & 1 &  1& -1&  0& 0  & 1 & 0 & 0 & 0 & -1 \\
 0 & 1 &  0 &  1& 0 & -1& 0  & 0 & 1 & 0 & 0 & -1\\
-1 & 1 &  1 &  0& 0 &  0& -1 & 0 & 0 & 1 & 0 & 0
\end{tabular}

$----------------------------------$

$V(11)$: $\deg = 12$, $h^{2,1} =7$, $rk =8$, $sq =14$, $dp =14$, $py =0$, 
${\rm vert}(\Delta^*) = 13$: 

\begin{tabular}{|c|c|c|c|c|c|c|c|c|c|c|c|c|}
1  & -1 & -1&  1& -1& -1& 1  & 1 &-1 & 1 & -1 &-1 &-1\\
-1 & 1  & 0 & -1&  1&  0&-1  & 0 & 1 &-1 & 0 & 1  & 0\\
-1 & 1 &  1 &  0& 0 &  0&-1  & -1& 1 & 0 & 1 & 0  & 0\\
-1 & 1 &  1 &  0& 1 &  1&  0 & -1& 0 &-1 & 0 & 0  & 0
\end{tabular}

$----------------------------------$

Type $(1,7)$:

$----------------------------------$

$V(12)$: $\deg = 14$, $h^{2,1} =5$, $rk =6$, $sq =9$, $dp =10$, $py =0$, 
${\rm vert}(\Delta^*) = 11$:

\begin{tabular}{|c|c|c|c|c|c|c|c|c|c|c|}
1 & -1 &-1&-1& 1&-1&-1  & 1 &-1 &-1&  1 \\
-1&  0 & 1& 1& 0& 1& 0  & 0 & 0 & 0& -1 \\
0 &  1 & 0& 1&-1& 0& 1  & 0 & 0 & 0& -1\\
-1 &  1 & 1& 0& 0& 0& 0  &-1 & 1 & 0&  0
\end{tabular}

$----------------------------------$

$V(13)$: $\deg = 14$, $h^{2,1} =5$, $rk =6$, $sq =10$, $dp =10$, $py =0$, 
${\rm vert}(\Delta^*) = 11$:

\begin{tabular}{|c|c|c|c|c|c|c|c|c|c|c|}
0 &  0 & 0& 0& 1& 0& 0  & 0 &-1 &-1& -1 \\
0 &  0 & 1& 0& 0&-1&-1  & 0 & 0 & 1& 0 \\
0 &  0 & 0& 1& 0&-1& 0  &-1 &0  & 0&  1\\
-1&  1 & 0& 0& 0& 1& 0  & 0 & 0 &-1& -1
\end{tabular}

$----------------------------------$

$V(14)$: 
 $\deg = 14$, $h^{2,1} =5$, $rk =7$, $sq =11$, $dp =11$, $py =0$, 
${\rm vert}(\Delta^*) = 12$:

\begin{tabular}{|c|c|c|c|c|c|c|c|c|c|c|c|}
1  & -1 & -1&  1& -1& -1& 1  & 1 &-1 &-1 & -1 &-1\\
-1 & 1  & 0 & -1& 1 &  0& -1 & 0 & 1 & 0 & 1  &0 \\
-1 & 1 &  1 &  0& 0 &  0& -1& -1 & 1 & 1 & 0 & 0\\
-1 & 1 &  1 &  0& 1 &  1& 0 & -1 & 0 & 0 & 0 & 0
\end{tabular}

$----------------------------------$

Type $(1,8)$

$----------------------------------$

$V(15)$:  $\deg = 16$, $h^{2,1} =3$, 
$rk =5$, $sq =7$, $dp =7$, $py =0$, ${\rm vert}(\Delta^*) = 10$:

\begin{tabular}{|c|c|c|c|c|c|c|c|c|c|}
 0 &-1& 0&-1 & 0& 0  & 1 & 0 & 0&  0 \\
-1 & 0& 1& 0 & 0& 0  & 0 &-1 & 0&  0 \\
-1 & 1& 0& 0 & 1& 0  & 0 & 0 & 0& -1\\
-1 & 1& 0& 0 & 0&-1  & 0 & 0 & 1&  0
\end{tabular}

$----------------------------------$

\newpage

$V(16)$:  $\deg = 16$, $h^{2,1} =3$,
 $rk =5$, $sq =5$, $dp =7$, $py =0$, ${\rm vert}(\Delta^*) = 10$:

\begin{tabular}{|c|c|c|c|c|c|c|c|c|c|}
-1 & 0&-1& 0 & 0& 1  & 0 &-1 & 0&  1 \\
-1 & 0& 1& 0 & 1&-1  & 0 & 0 &-1&  0 \\
 0 & 0& 1& 1 & 0&-1  &-1 & 0 & 0&  0\\
-1 & 1& 0& 0 & 0& 0  & 0 & 0 & 0&  0
\end{tabular}

$----------------------------------$

$V(17)$: $\deg = 16$, $h^{2,1} =3$,
 $rk =6$, $sq =8$, $dp =8$, $py =0$, ${\rm vert}(\Delta^*) = 11$:

\begin{tabular}{|c|c|c|c|c|c|c|c|c|c|c|}
-1 & 0 &-1& 1& 0& 0&-1  & 0 & 0 & 0& -1 \\
-1 & 0 & 0& 0&-1& 0&-1  & 0 & 1 & 0& 0 \\
0 &  0 &-1& 0& 0&-1&-1  & 1 & 0 & 0& 0\\
-1&  1 & 1& 0& 0& 0& 0  & 0 & 0 &-1& 0
\end{tabular}

$----------------------------------$

$V(18)$:  $\deg = 16$, $h^{2,1} =3$, 
$rk =7$, $sq =9$, $dp =9$, $py =0$, ${\rm vert}(\Delta^*) = 12$:

\begin{tabular}{|c|c|c|c|c|c|c|c|c|c|c|c|}
0 & -2 & 0 &  0& -1 &  1&  0&  0 & 0 & -1 &-1 &-1\\
0 & -1 & 1 &  0&  0 &  0&  0& -1 & 0 &  0 &-1 &-1 \\
0 & -1 & 0 &  1&  0 &  0& -1&  0 & 0 & -1 & 0 &-1\\
1 & -1 & 0 &  0&  0 &  0&  0&  0& -1 &-1  &-1 & 0
\end{tabular}

$----------------------------------$

Type $(1,9)$:

$----------------------------------$

$V(19)$: $\deg = 18$, $h^{2,1} =2$,  $rk =4$, 
$sq =4$, $dp =5$, $py =0$, ${\rm vert}(\Delta^*) =9$:

\begin{tabular}{|c|c|c|c|c|c|c|c|c|}
 0 & 0& 1& 0 & 0&-1  & 0 &-1 & 0  \\
-1 & 0& 0& 0 &-1& 1  & 1 & 0 & 0 \\
-1 & 0& 0& 1 & 0& 1  & 0 & 0 &-1 \\
-1 & 1& 0& 0 & 0& 0  & 0 & 0 & 0
\end{tabular}

$----------------------------------$

$V(20)$: $\deg = 18$, $h^{2,1} =2$, 
 $rk =5$, $sq =6$, $dp =6$, $py =0$, ${\rm vert}(\Delta^*) =10$:

\begin{tabular}{|c|c|c|c|c|c|c|c|c|c|}
-1 & 0& 0& 1 & 0&-1  &-1 & 0 &-1&  0 \\
 0 & 0& 1& 0 & 0&-1  & 0 &-1 &-1&  0 \\
-1 & 0& 0& 0 &-1& 1  & 0 & 0 & 0&  1\\
-1 & 1& 0& 0 & 0& 0  & 0 & 0 &-1&  0
\end{tabular}

$----------------------------------$

$V(21)$: $\deg = 18$, $h^{2,1} =2$, 
 $rk =5$, $sq =6$, $dp =6$, $py =0$, ${\rm vert}(\Delta^*) =10$:

\begin{tabular}{|c|c|c|c|c|c|c|c|c|c|}
 0 & 0& 1& 0 & 0&-1  & 0 & 0 &-1& -1 \\
-1 & 0& 0& 0 &-1& 0  & 1 & 0 & 0& -1 \\
-1 & 0& 0& 1 & 0& 0  & 0 &-1 &-1&  0\\
-1 & 1& 0& 0 & 0& 0  & 0 & 0 &-1& -1
\end{tabular}

$----------------------------------$

$V(22)$: $\deg = 18$, $h^{2,1} =2$,  
$rk =6$, $sq =7$, $dp =7$, $py =0$, ${\rm vert}(\Delta^*) =11$:

\begin{tabular}{|c|c|c|c|c|c|c|c|c|c|c|}
0 & -2 & 0& 0& -1& 1& 0  & 0 &-1 &-1& -1 \\
0 & -1 & 1& 0& 0 & 0& 0  &-1 &-1 &-1& 0 \\
0 & -1 & 0& 1& 0 & 0&-1  & 0 &0  &-1& -1\\
1 & -1 & 0& 0& 0 & 0& 0  & 0 &-1 & 0& -1
\end{tabular}

$----------------------------------$

Type $(1,11)$
 
$----------------------------------$

$V(23)$:   $\deg = 22$, $h^{2,1} =0$, $rk =9$, $sq =9$, $dp =9$, $py =1$, 
${\rm vert}(\Delta^*) =14$:

\begin{tabular}{|c|c|c|c|c|c|c|c|c|c|c|c|c|c|}
0 & -3 & -2&  0&  1 & -1& -1& -2 & -1&  0 &0  &-1 &-2 &-2 \\
-1&  0 & -1&  0&  0 & -1&  1&  0 & 1 &  1 &0 & -1 & 0 & 1\\
1 & -1 & 0 &  1&  0 &  1&  0&  0 & -1&  0 &0 & 0  &-1 & -1\\
1 & -1 & 0 &  0&  0 &  0& -1& -1 & 0 &  0 &1 & 0 & 0 & -1
\end{tabular}

$----------------------------------$

}

There $8$ types of smooth Fano $3$-folds with the Picard number $1$ 
that do not admit a toric 
conifold degeneration, but admit a conifold degeneration to a hypersurface 
in a toric variety. These types are the following:
$(2,2)$, $(2,3)$, $(1,4)$, $(1,5)$,  $(1,6)$, $(1,7)$,  $(1,8)$, $(1,9)$.

\section{The Golyshev correspondence}

In \cite{Go04} Golyshev has discovered a wonderful bijective correspondence
between $17$ types of  Fano $3$-folds with the Picard number $1$
and $17$ types of differential equations of  $D3$-type 
satisfying some modularity
conditions.

Let $Y$ be a Fano $3$-fold with Picard number
$1$. We put $H:= -K_Y$ and define the $4\times 4$-matix 

\[ A := \begin{pmatrix} a_{00} & a_{01} & a_{02} & a_{03} \\
1  & a_{11} & a_{12} & a_{13} \\
 0 & 1  & a_{22} & a_{23} \\
 0 & 0 & 1 & a_{33}
\end{pmatrix}
\]
where 
\[ a_{ij} = a_{3-j, 3-i} = \frac{j-i+1}{(-K_Y)^3}  \langle
H^{3-i}, H^j, \rangle_{j-i+1}. \]
and $\langle
H^{3-i}, H^j \rangle_{j-i+1}$ is  the number of maps
$f \, : \, \P^1 \to Y$ of degree $(j-i +1)$  such that
$f(0) \in L_{3-i}$, $f(\infty) \in L_j$ (here $L_k$ is the subvariety 
of codimension $k$ representing the cohomology class $H^k$). The matrix 
$A$ is called {\em counting matrix} of $Y$. 
One considers the equivalence class of such matrices up to adding a scalar 
matrix: $A \sim \lambda E_4 + A$. 
The matrix $A$ defines a $D3$-differential 
operator ${\mathcal D}_A$ expressed as polynomial 
in $D= t\frac{\partial}{\partial t}$: 

\[  {\mathcal D}_A : = 
D^3 - t(2D + 1) \left( (a_{00} + a_{11})D^2 + (a_{00}  +
a_{11}) D + a_{00} \right) + \]
\[ + t^2 (D+1) \left( ( a_{11}^2  + a_{00}^2  +
4 a_{11}a_{00} -a_{12} -2a_{01} ) D^2  + \right.\]
\[ \left. +
( 8 a_{11} a_{00}  -2a_{12} D - 4a_{01} + 2a_{11}^2)D +
6 a_{11} a_{00} - 4a_{01} \right) - \]
\[ -
t^3(2D+3)(D+2)(D+1) (
a_{00}^2 a_{11} + a_{11}^2 a_{00} - a_{12} a_{00} + a_{02}
- a_{11} a_{01} - a_{01} a_{00}) + \]
\[ +
t^4(D+3)(D+2)(D+1) ( - a_{00}^2 a_{12} + 2 a_{02} a_{00}
+ a_{00}^2 a_{11}^2 - a_{03} + a_{01}^2 - 2 a_{01}a_{11} a_{00})
\]
The $D3$-operators corresponding to equivalent matrices 
are called equivalent to each other.  The following result 
is due to Golyshev \cite{Go04}: 

\begin{theorem} 
There exists exactly $17$ types of $D3$-operators ${\mathcal D}_A$ such that 
the corresponding Picard-Fuchs equation ${\mathcal D}_A F(t) =0$ 
comes from the degree $d$ cyclic covering of a modular family over 
the universal elliptic curve with level $N$, so called $(d,N)$-modular 
family. The set of pairs of integers $(d,N)$ is exactly  
the set of all types smooth Fano $3$-folds $Y$ with the Picard number $1$.  
\end{theorem} 

It was proved case by case that the counting matrix for $17$  $(d,N)$-modular 
families is equivalent to the counting matrix of the Fano $3$-fold 
$Y$ (see \cite{Pr07a,Pr07b}). 

Our pupose was to find $D3$-operators from conifold degenerations 
$X \in {\mathcal F}(\Delta)$ as operators that annihilate 
the hypergeometric series $\Phi_0(t)$ defined in the previous 
section. 

Our result is the following:

\begin{theorem} 
For all $23$ $4$-dimensional reflexive polytopes 
$\Delta$ such that the Fano hypersurface $X \in {\mathcal F}(\Delta)$ 
has the Picard number $1$, the power series $\Phi_0(t)$ satisfies 
a $D3$-equation.  The matrices  of these $D3$-equations are equivalent 
to the matrices in the Golyshev's list (see  \cite{Go07}) and they 
coincide with them except for the case 
$(1,11)$ where the operator has the form
\[ D^3 -2 t D(1+D)(1 + 2D) - 8 t^2(1 + D) (12 + 22D + 
11D^2) -  \] 
\[ - 150 t^3 (1 + D) (2 + D) (3 + 2D) -
304 t^4 (1 + D) (2 + D) (3 + D) \]
and the corresponding matrix is
\[  \begin{pmatrix} 0 & 24 & 198 & 880 \\
1  & 2 & 44 & 198 \\
 0 & 1  & 2 & 24 \\
 0 & 0 & 1 & 0
\end{pmatrix}
\]
\end{theorem}

We expect that the multidimensional power series $\Phi(t_1, \ldots, t_r)$ 
can be used in a similar way to get Gromov-Witten invariant 
of Fano $3$-folds $Y$ with the Picard number $r \geq 2$ 
that admit a degeneration to a Fano hypersurface 
$X \in {\mathcal F}(\Delta)$. For this purpose we 
present the list of all  $143 =166-23$ polytopes in the case $r \geq 2$ 
and their invariants in the remaining sections.

\section{$B_2  =2$}

{\tiny

V(24):  $\deg = 12$,$\,$ $h^{1,2} = 9$, $rk =6$,
$sq =13$,
$dp =14$, $py =0$,  ${\rm vert}(\Delta^*) =12$:

\begin{tabular}{|c|c|c|c|c|c|c|c|c|c|c|c|}
-1&-1&-1&-1 & -1& 1 &-1 &-1 &-1 & 1 &-1 & 1 \\
 0& 0& 0& 1 &  0&-1 & 1 & 1 & 1 & 0 & 1 & 0\\
 0& 0& 1& 0 &  1& 0 & 1 & 1 & 0 &-1 & 0 &-1\\
 0& 1& 1& 0 &  0&-1 & 0 & 1 & 1 & 0 & 0 &-1
\end{tabular}

$----------------------------------$

V(25):  $\deg = 12$,$\,$ $h^{1,2} = 9$, $rk =6$,
$sq =14$,
$dp =14$, $py =0$,  ${\rm vert}(\Delta^*) =12$:

\begin{tabular}{|c|c|c|c|c|c|c|c|c|c|c|c|}
 1&-1& 1&-1 & -1&-1 &-1 &-1 &-1 &-1 & 1 & 1 \\
 0& 1& 0& 1 &  1& 0 & 1 & 0 & 0 & 0 &-1 &-1\\
 0& 1&-1& 0 &  1& 1 & 0 & 0 & 0 & 1 &-1 & 0\\
-1& 1&-1& 1 &  0& 1 & 0 & 2 & 1 & 0 & 0 &-1
\end{tabular}

$----------------------------------$

V(26):  $\deg = 14$,$\,$ $h^{1,2} = 5$, $rk =5$,
$sq =9$,
$dp =9$, $py =0$,  ${\rm vert}(\Delta^*) =11$:

\begin{tabular}{|c|c|c|c|c|c|c|c|c|c|c|}
 0& 0& 0&-1 &  0&-1 & 1 & 0 & 0 & 0 & 0 \\
 0& 0& 1& 0 &  1& 0 & 0 & 0 & 0 &-1 &-1 \\
 0& 1& 0& 0 & -1&-1 & 0 &-1 & 0 & 1 & 0 \\
-1& 0& 0& 0 & -1&-1 & 0 & 0 & 1 & 1 & 0
\end{tabular}

$----------------------------------$

V(27):  $\deg = 14$,$\,$ $h^{1,2} = 9$, $rk =4$,
$sq =9$,
$dp =12$, $py =0$,  ${\rm vert}(\Delta^*) =10$:

\begin{tabular}{|c|c|c|c|c|c|c|c|c|c|}
-1&-1&-1&-1 &  1& 1 &-1 & 1 &-1 & 1\\
 0& 0& 0& 0 &  0& 0 & 1 &-1 & 1 &-1\\
 0& 0& 1& 1 &  0& 0 & 0 &-1 & 0 &-1\\
 0& 1& 1& 0 &  0&-1 & 1 & 0 & 0 &-1
\end{tabular}

$----------------------------------$

V(28):  $\deg = 16$,$\,$ $h^{1,2} = 3$, $rk =4$,
$sq =5$,
$dp =6$, $py =0$,  ${\rm vert}(\Delta^*) =10$:

\begin{tabular}{|c|c|c|c|c|c|c|c|c|c|}
 0& 0& 0&-1 &  0& 0 & 1 & 0 & 1 &-1 \\
-1& 0& 0& 0 &  1& 0 & 0 &-1 &-1 & 1 \\
-1& 0&-1& 0 &  0& 1 & 0 & 0 &-1 & 1 \\
-1& 1& 0& 0 &  0& 0 & 0 & 0 & 0 & 0
\end{tabular}

$----------------------------------$

\newpage

V(29):  $\deg = 16$,$\,$ $h^{1,2} = 3$, $rk =4$,
$sq =6$,
$dp =6$, $py =0$,  ${\rm vert}(\Delta^*) =10$:

\begin{tabular}{|c|c|c|c|c|c|c|c|c|c|}
 0& 0& 1& 0 &  0& 1 & 0 &-1 & 0 &-1 \\
 0& 1& 0& 0 &  0&-1 &-1 & 0 & 0 & 1 \\
 0& 0& 0& 0 & -1&-1 & 0 & 0 & 1 & 1 \\
-1& 0& 0& 1 &  0& 0 & 0 & 0 & 0 & 0
\end{tabular}

$----------------------------------$

V(30):  $\deg = 16$,$\,$ $h^{1,2} = 5$, $rk =4$,
$sq =8$,
$dp =8$, $py =0$,  ${\rm vert}(\Delta^*) =10$:

\begin{tabular}{|c|c|c|c|c|c|c|c|c|c|}
-1& 0& 0& 0 &  0& 0 &-1 &-1 & 1 & 0 \\
-1&-1& 0& 1 &  1& 0 & 0 & 0 & 0 & 0 \\
 0& 0&-1& 0 &  1& 0 & 0 &-1 & 0 & 1 \\
-1& 0& 0& 0 &  1& 1 & 0 &-1 & 0 & 0
\end{tabular}

$----------------------------------$

V(31):  $\deg = 16$,$\,$ $h^{1,2} = 5$, $rk =4$,
$sq =7$,
$dp =8$, $py =0$,  ${\rm vert}(\Delta^*) =10$:

\begin{tabular}{|c|c|c|c|c|c|c|c|c|c|}
 0& 0&-1&-1 &  0& 0 &-1 & 0 & 1 & 0 \\
 1&-1&-1& 0 &  0& 0 & 0 & 1 & 0 & 0 \\
-1& 0& 1& 0 &  0&-1 &-1 & 0 & 0 & 1 \\
-1& 0& 0& 0 &  1& 0 &-1 & 0 & 0 & 0
\end{tabular}

$----------------------------------$

V(32):  $\deg = 16$,$\,$ $h^{1,2} = 5$, $rk =5$,
$sq =9$,
$dp =9$, $py =0$,  ${\rm vert}(\Delta^*) =11$:

\begin{tabular}{|c|c|c|c|c|c|c|c|c|c|c|}
 0& 0& 0&-1 &  0&-1 &-1 & 0 &-1 & 0 & 1 \\
 0& 0& 0& 0 & -1&-1 & 0 &-1 &-1 & 1 & 0 \\
-1& 1& 0& 0 & -1&-1 & 1 & 0 & 0 & 0 & 0 \\
 0& 0& 1& 0 &  1& 0 &-1 & 0 &-1 & 0 & 0
\end{tabular}

$----------------------------------$

V(33):  $\deg = 16$,$\,$ $h^{1,2} = 5$, $rk =5$,
$sq =9$,
$dp =9$, $py =0$,  ${\rm vert}(\Delta^*) =11$:

\begin{tabular}{|c|c|c|c|c|c|c|c|c|c|c|}
 0&-1&-1& 0 &  0&-1 &-1 & 0 & 0 & 0 & 1 \\
 0& 0&-1& 0 & -1& 0 &-1 &-1 & 1 & 0 & 0 \\
 0& 0& 1&-1 &  0&-1 & 0 &-1 & 0 & 1 & 0 \\
 1& 0& 0& 0 &  0&-1 &-1 &-1 & 0 & 0 & 0
\end{tabular}

$----------------------------------$

V(34):  $\deg = 18$,$\,$ $h^{1,2} = 5$, $rk =3$,
$sq =5$,
$dp =7$, $py =0$,  ${\rm vert}(\Delta^*) =9$:

\begin{tabular}{|c|c|c|c|c|c|c|c|c|}
-1& 1& 1& 0 &  0&-1& 0 & 0 & 0\\
 0& 0& 0& 0 & -1&-1&-1 & 1 & 0\\
 0& 0&-1& 0 &  0& 0&-1 & 0 & 1\\
 0& 0&-1& 1 &  0& 1&-1 & 0 & 0
\end{tabular}

$----------------------------------$

V(35):  $\deg = 18$,$\,$ $h^{1,2} = 5$, $rk =4$,
$sq =7$,
$dp =8$, $py =0$,  ${\rm vert}(\Delta^*) =10$:

\begin{tabular}{|c|c|c|c|c|c|c|c|c|c|}
-1&-1& 1& 0 &  0&-1 &-1 & 0 & 0 & 0 \\
 0& 0& 0& 0 & -1&-1 &-1 &-1 & 1 & 0 \\
-1& 0& 0& 0 &  0& 0 &-1 &-1 & 0 & 1 \\
-1& 0& 0& 1 &  0& 1 & 0 &-1 & 0 & 0
\end{tabular}

$----------------------------------$

V(36):  $\deg = 20$,$\,$ $h^{1,2} = 1$, $rk =3$,
$sq =3$,
$dp =3$, $py =0$,  ${\rm vert}(\Delta^*) =9$:

\begin{tabular}{|c|c|c|c|c|c|c|c|c|}
 0& 0& 0& 0 &  1&-1&-1 & 0 & 0\\
 0& 0& 0&-1 &  0& 0&-1 & 0 & 1\\
 0& 0& 1& 0 &  0& 0&-1 &-1 & 0\\
-1& 1& 0& 0 &  0& 0& 0 & 0 & 0
\end{tabular}

$----------------------------------$

V(37):  $\deg = 20$,$\,$ $h^{1,2} = 2$, $rk =3$,
$sq =4$,
$dp =4$, $py =0$,  ${\rm vert}(\Delta^*) =9$:

\begin{tabular}{|c|c|c|c|c|c|c|c|c|}
-1& 0& 0& 0 &  0&-1&-1 & 0 & 1\\
 0&-1& 0& 0 &  0& 0&-1 & 1 & 0\\
-1& 0&-1& 0 &  1& 0& 0 & 0 & 0\\
-1& 0& 0& 1 &  0& 0&-1 & 0 & 0
\end{tabular}

$----------------------------------$

V(38):  $\deg = 20$,$\,$ $h^{1,2} = 2$, $rk =3$,
$sq =3$,
$dp =4$, $py =0$,  ${\rm vert}(\Delta^*) =9$:

\begin{tabular}{|c|c|c|c|c|c|c|c|c|}
-1& 1& 0& 0 &  0&-1&-1 & 0 & 0\\
 0& 0& 0& 1 &  0& 0&-1 & 0 &-1\\
-1& 0& 0& 0 &  1& 0& 1 &-1 & 0\\
-1& 0& 1& 0 &  0& 0& 0 & 0 & 0
\end{tabular}

$----------------------------------$

\newpage

V(39):  $\deg = 20$,$\,$ $h^{1,2} = 2$, $rk =4$,
$sq =5$,
$dp =5$, $py =0$,  ${\rm vert}(\Delta^*) =10$:

\begin{tabular}{|c|c|c|c|c|c|c|c|c|c|}
-1& 0& 0&-1 &  0&-1 &-2 & 0 & 1 & 0 \\
 0& 1& 0&-1 & -1& 0 &-1 & 0 & 0 & 0 \\
-1& 0& 0&-1 &  0& 0 &-1 &-1 & 0 & 1 \\
-1& 0& 1& 0 &  0& 0 &-1 & 0 & 0 & 0
\end{tabular}

$----------------------------------$

V(40):  $\deg = 20$,$\,$ $h^{1,2} = 3$, $rk =2$,
$sq =2$,
$dp =4$, $py =0$,  ${\rm vert}(\Delta^*) =8$:

\begin{tabular}{|c|c|c|c|c|c|c|c|}
 0& 0& 1& 0 &  0&-1 &-1 & 1\\
-1& 0& 0& 0 &  1& 0 &-1 & 1\\
 0& 0& 0& 1 &  0& 0 &-1 & 0\\
 0& 1& 0& 0 &  0& 0 & 0 &-1
\end{tabular}

$----------------------------------$

V(41):  $\deg = 20$,$\,$ $h^{1,2} = 3$, $rk =3$,
$sq =5$,
$dp =5$, $py =0$,  ${\rm vert}(\Delta^*) =9$:

\begin{tabular}{|c|c|c|c|c|c|c|c|c|}
-1& 0&-1&-1 &  0& 0& 0 &-1 & 1\\
 0& 1& 0& 1 & -1& 0& 0 &-1 & 0\\
-1& 0& 0& 0 &  0& 1& 0 &-1 & 0\\
-1& 0& 0&-1 &  0& 0& 1 & 0 & 0
\end{tabular}

$----------------------------------$

V(42):  $\deg = 20$,$\,$ $h^{1,2} = 3$, $rk =3$,
$sq =4$,
$dp =5$, $py =0$,  ${\rm vert}(\Delta^*) =9$:

\begin{tabular}{|c|c|c|c|c|c|c|c|c|}
-1& 0&-1&-1 &  0& 0& 1 & 0 & 0\\
 0& 0& 0&-1 &  1& 0& 0 &-1 &-1\\
-1& 0& 0&-1 &  0& 1& 0 & 0 &-1\\
-1& 1& 0& 0 &  0& 0& 0 & 0 &-1
\end{tabular}

$----------------------------------$

V(43):  $\deg = 20$,$\,$ $h^{1,2} = 3$, $rk =4$,
$sq =6$,
$dp =6$, $py =0$,  ${\rm vert}(\Delta^*) =10$:

\begin{tabular}{|c|c|c|c|c|c|c|c|c|c|}
-1& 0& 0&-1 &  0& 1 &-1 & 0 & 0 &-1 \\
-1&-1& 0& 0 &  1& 0 & 0 & 0 &-1 &-1 \\
 0& 0& 1& 1 &  0& 0 & 0 & 0 &-1 &-1 \\
-1& 0& 0&-1 &  0& 0 & 0 & 1 & 1 & 0
\end{tabular}

$----------------------------------$

V(44):  $\deg = 20$,$\,$ $h^{1,2} = 3$, $rk =4$,
$sq =6$,
$dp =6$, $py =0$,  ${\rm vert}(\Delta^*) =10$:

\begin{tabular}{|c|c|c|c|c|c|c|c|c|c|}
-1& 0&-1& 0 & -1&-2 &-1 & 0 & 0 & 1 \\
-1& 1&-1& 0 &  0&-1 & 0 &-1 & 0 & 0 \\
 0& 0&-1& 0 &  0&-1 &-1 & 0 & 1 & 0 \\
-1& 0& 0& 1 &  0&-1 &-1 & 0 & 0 & 0
\end{tabular}

$----------------------------------$

V(45):  $\deg = 20$,$\,$ $h^{1,2} = 3$, $rk =9$,
$sq =12$,
$dp =12$, $py =1$,  ${\rm vert}(\Delta^*) =15$:

\begin{tabular}{|c|c|c|c|c|c|c|c|c|c|c|c|c|c|c|}
-1&-1&-1&-1 & -1&-1 &-1 &-1 &-1 &-1 &-1 &-1& -1&-1& 1 \\
 1& 1& 0& 0 &  1& 2 & 2 & 2 & 2 & 1 & 0 & 0&  1& 1&-1\\
 1& 0& 1& 0 &  0& 0 & 1 & 2 & 1 & 1 & 1 & 2&  2& 2&-1\\
 2& 1& 1& 0 &  0& 1 & 2 & 2 & 1 & 0 & 0 & 1&  1& 2&-1
\end{tabular}

$----------------------------------$

V(46):  $\deg = 22$,$\,$ $h^{1,2} = 2$, $rk =2$,
$sq =2$,
$dp =3$, $py =0$,  ${\rm vert}(\Delta^*) =8$:

\begin{tabular}{|c|c|c|c|c|c|c|c|}
 1& 0& 1& 0 &  0& 0 &-1 & 0\\
 0& 0& 0& 1 &  1& 0 & 0 &-1\\
-1& 0& 0& 0 & -1& 1 & 0 & 0\\
-1& 1& 0& 0 & -1& 0 & 0 & 0
\end{tabular}

$----------------------------------$

V(47):  $\deg = 22$,$\,$ $h^{1,2} = 2$, $rk =3$,
$sq =4$,
$dp =4$, $py =0$,  ${\rm vert}(\Delta^*) =9$:

\begin{tabular}{|c|c|c|c|c|c|c|c|c|}
-1& 0& 1& 0 &  0& 0&-1 &-1 &-2\\
 0& 1& 0& 0 &  0&-1& 0 &-1 &-1\\
-1& 0& 0& 1 &  0& 0& 0 &-1 &-1\\
-1& 0& 0& 0 &  1& 0& 0 & 0 &-1
\end{tabular}

$----------------------------------$

V(48):  $\deg = 22$,$\,$ $h^{1,2} = 2$, $rk =4$,
$sq =5$,
$dp =5$, $py =0$,  ${\rm vert}(\Delta^*) =10$:

\begin{tabular}{|c|c|c|c|c|c|c|c|c|c|}
-1& 0& 0& 1 &  0&-1 &-1 &-2 & 0 &-1 \\
-1& 0& 1& 0 &  0& 0 &-1 &-1 &-1 &-2 \\
 0& 1& 0& 0 &  0& 0 &-1 &-1 & 0 &-1 \\
-1& 0& 0& 0 &  1& 0 & 0 &-1 & 0 &-1
\end{tabular}

$----------------------------------$

\newpage

V(49):  $\deg = 22$,$\,$ $h^{1,2} = 4$, $rk =2$,
$sq =3$,
$dp =5$, $py =0$,  ${\rm vert}(\Delta^*) =8$:

\begin{tabular}{|c|c|c|c|c|c|c|c|}
-1& 0& 0& 0 &  1& 1 & 1 & 0\\
 0& 0& 1& 0 & -1& 0 &-1 &-1\\
 0& 0& 0& 1 & -1& 0 & 0 &-1\\
 0& 1& 0& 0 &  0& 0 &-1 &-1
\end{tabular}

$----------------------------------$

V(50):  $\deg = 22$,$\,$ $h^{1,2} = 4$, $rk =3$,
$sq =6$,
$dp =6$, $py =0$,  ${\rm vert}(\Delta^*) =9$:

\begin{tabular}{|c|c|c|c|c|c|c|c|c|}
-2&-1& 0&-1 & -1&-1& 0 & 0 & 1\\
-1& 0& 0&-1 &  0&-1& 1 & 0 & 0\\
-1& 0& 0&-1 & -1& 0& 0 & 1 & 0\\
-1& 0& 1& 0 & -1&-1& 0 & 0 & 0
\end{tabular}

$----------------------------------$

V(51):  $\deg = 24$,$\,$ $h^{1,2} = 1$, $rk =2$,
$sq =2$,
$dp =2$, $py =0$,  ${\rm vert}(\Delta^*) =8$:

\begin{tabular}{|c|c|c|c|c|c|c|c|}
-1& 0&-1&-1 &  0& 0 & 0 & 1\\
 1& 0& 0&-1 & -1& 1 & 0 & 0\\
 0& 0& 0&-1 &  0& 0 & 1 & 0\\
-1& 1& 0& 0 &  0& 0 & 0 & 0
\end{tabular}

$----------------------------------$

V(52):  $\deg = 24$,$\,$ $h^{1,2} = 1$, $rk =2$,
$sq =2$,
$dp =2$, $py =0$,  ${\rm vert}(\Delta^*) =8$:

\begin{tabular}{|c|c|c|c|c|c|c|c|}
 0& 0& 0&-1 & -1&-1 & 1 & 0\\
-1& 0& 1& 1 &  0& 0 & 0 & 0\\
 0& 0& 0&-1 &  0&-1 & 0 & 1\\
 0& 1& 0& 0 &  0&-1 & 0 & 0
\end{tabular}

$----------------------------------$

V(53):  $\deg = 24$,$\,$ $h^{1,2} = 1$, $rk =3$,
$sq =3$,
$dp =3$, $py =0$,  ${\rm vert}(\Delta^*) =9$:

\begin{tabular}{|c|c|c|c|c|c|c|c|c|}
-1& 0& 0&-1 & -1&-2& 0 & 0 & 1\\
-1& 0& 0&-1 &  0&-1&-1 & 1 & 0\\
 0& 1& 0&-1 &  0&-1& 0 & 0 & 0\\
-1& 0& 1& 0 &  0&-1& 0 & 0 & 0
\end{tabular}

$----------------------------------$

V(54):  $\deg = 24$,$\,$ $h^{1,2} = 1$, $rk =7$,
$sq =8$,
$dp =8$, $py =1$,  ${\rm vert}(\Delta^*) =13$:

\begin{tabular}{|c|c|c|c|c|c|c|c|c|c|c|c|c|}
-1& 0& 1& 0 & -1&-3 &-2 &-1 &-2 & 0 &-1 &-2& -2 \\
-1& 0& 0& 0 &  1& 0 & 0 &-1 &-1 & 1 & 1 & 1&  0\\
 1& 1& 0& 0 & -1&-1 &-1 & 0 & 0 & 0 & 0 &-1&  0\\
 0& 0& 0& 1 &  0&-1 & 0 & 1 & 0 & 0 &-1 &-1& -1
\end{tabular}

$----------------------------------$

V(55):  $\deg = 24$,$\,$ $h^{1,2} = 2$, $rk =2$,
$sq =2$,
$dp =3$, $py =0$,  ${\rm vert}(\Delta^*) =8$:

\begin{tabular}{|c|c|c|c|c|c|c|c|}
 0& 1& 0& 0 &  0&-1 & 0 &-1\\
-1& 0& 0& 1 &  0& 0 &-1 &-2\\
-1& 0& 0& 0 &  1& 0 & 0 &-1\\
-1& 0& 1& 0 &  0& 0 & 0 &-1
\end{tabular}

$----------------------------------$

V(56):  $\deg = 24$,$\,$ $h^{1,2} = 2$, $rk =3$,
$sq =3$,
$dp =4$, $py =0$,  ${\rm vert}(\Delta^*) =9$:

\begin{tabular}{|c|c|c|c|c|c|c|c|c|}
 0&-1& 0&-1 & -2&-1& 0 & 0 & 1\\
 0&-1&-1& 0 & -1&-2& 1 & 0 & 0\\
 0&-1& 0& 0 & -1&-1& 0 & 1 & 0\\
 1& 0& 0& 0 & -1&-1& 0 & 0 & 0
\end{tabular}

$----------------------------------$

V(57):  $\deg = 26$,$\,$ $h^{1,2} = 0$, $rk =6$,
$sq =6$,
$dp =6$, $py =1$,  ${\rm vert}(\Delta^*) =12$:

\begin{tabular}{|c|c|c|c|c|c|c|c|c|c|c|c|}
 1& 0&-1& 0 & -2&-2 & 0 & 1 &-1 &-2 &-2 & 0 \\
 0&-1&-1& 0 &  0& 1 & 1 & 0 &-1 &-1 & 0 & 0\\
 0& 1& 0& 1 & -1&-1 & 0 & 1 & 1 & 0 & 0 & 0\\
 0& 1& 1& 0 &  0&-1 & 0 & 1 & 0 & 0 &-1 & 1
\end{tabular}

$----------------------------------$

V(58):  $\deg = 26$,$\,$ $h^{1,2} = 0$, $rk =6$,
$sq =6$,
$dp =6$, $py =1$,  ${\rm vert}(\Delta^*) =12$:

\begin{tabular}{|c|c|c|c|c|c|c|c|c|c|c|c|}
-1& 0& 0& 1 & -1& 0 & 0 &-1 &-3 &-2 &-2 &-2 \\
 1& 1& 1& 0 & -1& 0 & 0 & 0 &-1 & 0 &-1 & 0\\
-1&-1& 0& 0 &  1& 1 & 0 &-1 & 0 & 0 & 0 &-1\\
 0& 1& 0& 0 &  0& 0 & 1 & 1 &-1 &-1 & 0 & 0
\end{tabular}

$----------------------------------$
\newpage

V(59):  $\deg = 26$,$\,$ $h^{1,2} = 2$, $rk =1$,
$sq =1$,
$dp =2$, $py =0$,  ${\rm vert}(\Delta^*) =7$:

\begin{tabular}{|c|c|c|c|c|c|c|}
-1& 0& 0& 0& 1 &  1& 0 \\
 0& 0& 1& 0& 0 & -1&-1\\
 0& 0& 0& 1& 0 & -1&-1\\
 0& 1& 0& 0& 0 &  0&-1
\end{tabular}

$----------------------------------$

V(60):  $\deg = 26$,$\,$ $h^{1,2} = 2$, $rk =2$,
$sq =3$,
$dp =3$, $py =0$,  ${\rm vert}(\Delta^*) =8$:

\begin{tabular}{|c|c|c|c|c|c|c|c|}
-1&-1& 0& 0 &  1& 0 &-1 &-1\\
-1& 0& 0& 1 &  0& 0 & 0 &-1\\
 0& 1& 1& 0 &  0& 0 & 0 &-1\\
-1&-1& 0& 0 &  0& 1 & 0 & 0
\end{tabular}

$----------------------------------$

V(61):  $\deg = 26$,$\,$ $h^{1,2} = 2$, $rk =6$,
$sq =8$,
$dp =8$, $py =1$,  ${\rm vert}(\Delta^*) =12$:

\begin{tabular}{|c|c|c|c|c|c|c|c|c|c|c|c|}
-1&-1&-1&-1 & -1&-1 &-1 &-1 &-1 &-1 &-1 & 1 \\
 0& 1& 0& 0 &  2& 1 & 1 & 1 & 0 & 2 & 1 &-1\\
 0& 0& 1& 0 &  1& 2 & 1 & 2 & 1 & 1 & 0 &-1\\
 1& 2& 1& 0 &  1& 0 & 0 & 1 & 0 & 2 & 1 &-1
\end{tabular}

$----------------------------------$

V(62):  $\deg = 28$,$\,$ $h^{1,2} = 0$, $rk =5$,
$sq =5$,
$dp =5$, $py =1$,  ${\rm vert}(\Delta^*) =11$:

\begin{tabular}{|c|c|c|c|c|c|c|c|c|c|c|}
 1&-1& 0&-2 & -2& 0 & 1 & 0 &-2 & 0 &-2 \\
 0& 1& 0&-1 &  0&-1 & 0 & 1 & 0 & 0 & 1 \\
 0& 0& 0& 0 & -1& 1 & 1 & 0 & 0 & 1 &-1 \\
 0&-1& 1& 0 &  0& 1 & 1 & 0 &-1 & 0 &-1
\end{tabular}

$----------------------------------$

V(63):  $\deg = 28$,$\,$ $h^{1,2} = 0$, $rk =5$,
$sq =5$,
$dp =5$, $py =1$,  ${\rm vert}(\Delta^*) =11$:

\begin{tabular}{|c|c|c|c|c|c|c|c|c|c|c|}
-1& 0& 1& 0 & -2&-1 &-3 &-1 &-2 & 0 &-2 \\
 1& 1& 0& 0 & -1&-1 &-1 & 0 & 0 & 0 & 0 \\
-1& 0& 0& 0 &  0& 1 & 0 &-1 &-1 & 1 & 0 \\
 0& 0& 0& 1 &  0& 0 &-1 & 1 & 0 & 0 &-1
\end{tabular}

$----------------------------------$

V(64):  $\deg = 28$,$\,$ $h^{1,2} = 0$, $rk =6$,
$sq =6$,
$dp =6$, $py =1$,  ${\rm vert}(\Delta^*) =12$:

\begin{tabular}{|c|c|c|c|c|c|c|c|c|c|c|c|}
-4& 0&-1& 0 & -2&-2 & 0 & 1 &-1 &-3 &-2 &-3 \\
-1& 0& 0& 1 &  0&-1 & 0 & 0 & 0 &-1 & 0 &-1\\
-1& 1& 1& 0 &  0& 0 & 0 & 0 &-1 &-1 &-1 & 0\\
-1& 0&-1& 0 & -1& 0 & 1 & 0 & 1 & 0 & 0 &-1
\end{tabular}

$----------------------------------$

V(65):  $\deg = 30$,$\,$ $h^{1,2} = 0$, $rk =4$,
$sq =4$,
$dp =4$, $py =1$,  ${\rm vert}(\Delta^*) =10$:

\begin{tabular}{|c|c|c|c|c|c|c|c|c|c|}
 1& 1& 0& 0 & -2& 0 &-2 &-1 &-2 &-2 \\
 0& 1& 0& 1 & -1& 0 &-1 & 1 & 0 & 0 \\
 0& 0& 1& 0 &  0& 0 & 1 &-1 &-1 & 0 \\
 0& 1& 0& 0 &  0& 1 &-1 & 0 & 0 &-1
\end{tabular}

$----------------------------------$

V(66):  $\deg = 30$,$\,$ $h^{1,2} = 0$, $rk =4$,
$sq =4$,
$dp =4$, $py =1$,  ${\rm vert}(\Delta^*) =10$:

\begin{tabular}{|c|c|c|c|c|c|c|c|c|c|}
-1& 0&-2& 0 &  1& 0 & 0 &-1 &-1 &-2 \\
 1& 1&-1& 1 &  0& 0 & 0 &-1 & 0 & 0 \\
 0& 0& 0& 1 &  0& 0 & 1 & 1 &-1 &-1 \\
-1& 0& 0&-1 &  0& 1 & 0 & 0 & 1 & 0
\end{tabular}

$----------------------------------$

V(67):  $\deg = 30$,$\,$ $h^{1,2} = 0$, $rk =5$,
$sq =5$,
$dp =5$, $py =1$,  ${\rm vert}(\Delta^*) =11$:

\begin{tabular}{|c|c|c|c|c|c|c|c|c|c|c|}
-2& 0&-1& 1 & -1& 0 &-2 &-2 & 0 &-1 &-2 \\
-1& 0&-1& 0 & -1& 0 &-1 & 0 & 1 & 0 & 0 \\
 0& 1& 1& 0 &  0& 0 &-1 &-1 & 0 & 1 & 0 \\
 0& 0& 0& 0 &  1& 1 & 1 & 0 & 0 &-1 &-1
\end{tabular}

$----------------------------------$

V(68):  $\deg = 30$,$\,$ $h^{1,2} = 1$, $rk =1$,
$sq =1$,
$dp =1$, $py =0$,  ${\rm vert}(\Delta^*) =7$:

\begin{tabular}{|c|c|c|c|c|c|c|}
 1& 0& 1& 0& 0 & -1& 1 \\
-1& 0& 0& 1& 0 &  0&-1\\
-1& 0& 0& 0& 1 &  0& 0\\
 0& 1& 0& 0& 0 &  0&-1
\end{tabular}

$----------------------------------$

\newpage 

V(69):  $\deg = 30$,$\,$ $h^{1,2} = 1$, $rk =4$,
$sq =5$,
$dp =5$, $py =1$,  ${\rm vert}(\Delta^*) =10$:

\begin{tabular}{|c|c|c|c|c|c|c|c|c|c|}
 1& 0& 0&-2 & -2& 1 & 0 & 0 &-2 &-2 \\
 0& 0& 0&-1 &  0& 0 &-1 & 1 & 0 & 1 \\
 0& 1& 0& 0 & -1& 1 & 1 & 0 & 0 &-1 \\
 0& 0& 1& 0 &  0& 1 & 1 & 0 &-1 &-1
\end{tabular}

$----------------------------------$

V(70):  $\deg = 32$,$\,$ $h^{1,2} = 1$, $rk =0$, $sq =0$,
$dp =0$, $py =0$, ${\rm vert}(\Delta^*) =6$:

\begin{tabular}{|c|c|c|c|c|c|}
-1& 0& 0& 0& 1 &  0 \\
 0& 0& 1& 0& 0 & -1\\
 0& 0& 0& 1& 0 & -1 \\
 0& 1& 0& 0& 0 & -1
\end{tabular}

$----------------------------------$

V(71):  $\deg = 32$,$\,$ $h^{1,2} = 1$, $rk =1$,
$sq =1$,
$dp =1$, $py =0$,  ${\rm vert}(\Delta^*) =7$:

\begin{tabular}{|c|c|c|c|c|c|c|}
-2& 0&-1&-1& 0 &  0& 1 \\
-1& 0& 0&-1& 1 &  0& 0\\
-1& 0& 0&-1& 0 &  1& 0\\
-1& 1& 0& 0& 0 &  0& 0
\end{tabular}

$----------------------------------$

V(72):  $\deg = 32$,$\,$ $h^{1,2} = 1$, $rk =3$,
$sq =4$,
$dp =4$, $py =1$,  ${\rm vert}(\Delta^*) =9$:

\begin{tabular}{|c|c|c|c|c|c|c|c|c|}
-2& 0& 0&-1 &  0&-2&-3 & 0 & 1\\
 0& 1& 0& 0 &  0&-1&-1 & 1 & 0\\
-1&-1& 1& 1 &  0& 0& 0 & 0 & 0\\
 0& 1& 0&-1 &  1& 0&-1 & 0 & 0
\end{tabular}

$----------------------------------$

V(73):  $\deg = 32$,$\,$ $h^{1,2} = 1$, $rk =4$,
$sq =5$,
$dp =5$, $py =1$,  ${\rm vert}(\Delta^*) =10$:

\begin{tabular}{|c|c|c|c|c|c|c|c|c|c|}
-3& 0&-1&-2 &  0& 0 &-2 & 0 &-3 & 1 \\
-1&-1&-1&-1 &  0& 1 & 0 & 0 & 0 & 0 \\
 0& 1& 1& 0 &  1& 0 &-1 & 0 &-1 & 0 \\
-1& 1& 0& 0 &  0& 0 & 0 & 1 &-1 & 0
\end{tabular}

$----------------------------------$

V(74):  $\deg = 34$,$\,$ $h^{1,2} = 0$, $rk =3$,
$sq =3$,
$dp =3$, $py =1$,  ${\rm vert}(\Delta^*) =9$:

\begin{tabular}{|c|c|c|c|c|c|c|c|c|}
 1& 0& 1&-2 &  0& 0&-2 &-4 &-2\\
 0& 0& 1& 0 &  1& 0&-1 &-1 & 0\\
 0& 0& 1&-1 &  0& 1& 0 &-1 & 0\\
 0& 1& 0& 0 &  0& 0& 0 &-1 &-1
\end{tabular}

$----------------------------------$

V(75):  $\deg = 34$,$\,$ $h^{1,2} = 0$, $rk =4$,
$sq =4$,
$dp =4$, $py =1$,  ${\rm vert}(\Delta^*) =10$:

\begin{tabular}{|c|c|c|c|c|c|c|c|c|c|}
-3&-2&-4& 0 & -1& 0 & 1 &-2 &-2 & 0 \\
-1& 0&-1& 0 & -1& 1 & 0 & 0 &-1 & 0 \\
 0& 0&-1& 1 &  1& 0 & 0 &-1 & 0 & 0 \\
-1&-1&-1& 0 &  0& 0 & 0 & 0 & 0 & 1
\end{tabular}

$----------------------------------$

V(76):  $\deg = 34$,$\,$ $h^{1,2} = 0$, $rk =5$,
$sq =5$,
$dp =5$, $py =1$,  ${\rm vert}(\Delta^*) =11$:

\begin{tabular}{|c|c|c|c|c|c|c|c|c|c|c|}
-5& 0&-2&-3 & -2&-4 & 0 & 0 &-3 &-2 & 1 \\
-1& 0&-1& 0 &  0&-1 & 1 & 0 &-1 & 0 & 0 \\
-1& 1& 0&-1 & -1&-1 & 0 & 0 & 0 & 0 & 0 \\
-2& 0& 0&-1 &  0&-1 & 0 & 1 &-1 &-1 & 0
\end{tabular}

$----------------------------------$

V(77):  $\deg = 38$,$\,$ $h^{1,2} = 0$, $rk =2$,
$sq =2$,
$dp =2$, $py =1$,  ${\rm vert}(\Delta^*) =8$:

\begin{tabular}{|c|c|c|c|c|c|c|c|}
-3& 0& 1& 0 & -1& -2 & 0 &-3\\
 0& 1& 0& 0 &  1& -1 & 0 &-1\\
-1& 0& 0& 1 &  0&  0 & 0 &-1\\
-1& 0& 0& 0 & -1&  0 & 1 & 0
\end{tabular}

$----------------------------------$

V(78):  $\deg = 38$,$\,$ $h^{1,2} = 0$, $rk =3$,
$sq =3$,
$dp =3$, $py =1$,  ${\rm vert}(\Delta^*) =9$:

\begin{tabular}{|c|c|c|c|c|c|c|c|c|}
-2&-1&-1&-3 &  0& 0& 1 & 0 &-2\\
 1& 1& 1& 0 &  1& 0& 0 & 0 &-1\\
-1&-1& 0&-1 &  0& 0& 0 & 1 & 0\\
-1& 0&-1&-1 &  0& 1& 0 & 0 & 0
\end{tabular}

$----------------------------------$
\newpage 

V(79):  $\deg = 38$,$\,$ $h^{1,2} = 0$, $rk =3$,
$sq =3$,
$dp =3$, $py =1$,  ${\rm vert}(\Delta^*) =9$:

\begin{tabular}{|c|c|c|c|c|c|c|c|c|}
 1& 0& 1&-3 &  0&-2 & 0 &-4&-3\\
 0& 0& 0&-1 &  0&-1 & 1 &-1&-1\\
 0& 0& 1&-1 &  1& 0 & 0 &-1& 0\\
 0& 1& 1& 0 &  0& 0 & 0 &-1&-1
\end{tabular}

$----------------------------------$

V(80):  $\deg = 40$,$\,$ $h^{1,2} = 0$, $rk =2$,
$sq =2$,
$dp =2$, $py =1$,  ${\rm vert}(\Delta^*) =8$:

\begin{tabular}{|c|c|c|c|c|c|c|c|}
-1& 0& 1& 0 & -2& 0 &-2 &-1\\
 1& 1& 0& 0 & -1& 0 & 0 & 0\\
 0& 0& 0& 1 &  0& 0 &-1 & 1\\
-1& 0& 0& 0 &  0& 1 & 0 &-1
\end{tabular}

$----------------------------------$

V(81):  $\deg = 40$,$\,$ $h^{1,2} = 0$, $rk =3$,
$sq =3$,
$dp =3$, $py =1$,  ${\rm vert}(\Delta^*) =9$:

\begin{tabular}{|c|c|c|c|c|c|c|c|c|}
-3& 0& 0&-2 & -2&-1&-3 & 0 & 1\\
-1& 0& 0&-1 &  0&-1&-2 & 1 & 0\\
 0& 1& 0& 0 & -1& 1& 1 & 0 & 0\\
-1& 0& 1& 0 &  0& 0&-1 & 0 & 0
\end{tabular}

$----------------------------------$

V(82):  $\deg = 40$,$\,$ $h^{1,2} = 1$, $rk =0$, $sq =0$,
$dp =0$, $py =0$,  ${\rm vert}(\Delta^*) =6$:

\begin{tabular}{|c|c|c|c|c|c|}
 2& 1& 0& 0& 0 & -1 \\
-1& 0& 0& 1& 0 &  0\\
-1& 0& 0& 0& 1 &  0 \\
-1& 0& 1& 0& 0 &  0
\end{tabular}

$----------------------------------$

V(83):  $\deg = 40$,$\,$ $h^{1,2} = 1$, $rk =2$,
$sq =3$,
$dp =3$, $py =1$,  ${\rm vert}(\Delta^*) =8$:

\begin{tabular}{|c|c|c|c|c|c|c|c|}
-1&-2& 1& 0 &  0& 1 & 1 &-2\\
 0& 0& 0& 0 &  1& 1 & 1 &-1\\
 0&-1& 0& 1 &  0& 1 & 1 & 0\\
 1&-1& 0& 0 &  0& 1 & 2 &-1
\end{tabular}

$----------------------------------$

V(84):  $\deg = 46$,$\,$ $h^{1,2} = 0$, $rk =1$,
$sq =1$,
$dp =1$, $py =1$,  ${\rm vert}(\Delta^*) =7$:

\begin{tabular}{|c|c|c|c|c|c|c|}
-3& 0& 0&-1& 1 &  0&-2 \\
 0& 1& 0& 1& 0 &  0&-1\\
-1& 0& 0&-1& 0 &  1& 0\\
-1& 0& 1& 0& 0 &  0& 0
\end{tabular}

$----------------------------------$

V(85):  $\deg = 46$,$\,$ $h^{1,2} = 0$, $rk =1$,
$sq =1$,
$dp =1$, $py =1$,  ${\rm vert}(\Delta^*) =7$:

\begin{tabular}{|c|c|c|c|c|c|c|}
-3&-1& 0& 0& 1 &  0&-3 \\
-1& 0& 0& 1& 0 &  0&-1\\
 0& 1& 1& 0& 0 &  0&-1\\
-1&-1& 0& 0& 0 &  1& 0
\end{tabular}

$----------------------------------$

V(86):  $\deg = 46$,$\,$ $h^{1,2} = 0$, $rk =2$,
$sq =2$,
$dp =2$, $py =1$,  ${\rm vert}(\Delta^*) =8$:

\begin{tabular}{|c|c|c|c|c|c|c|c|}
-3& 0&-2& 0 &  1& -2 & 0 &-3\\
-1& 0&-1& 1 &  0& -1 & 0 &-1\\
 0& 1& 1& 0 &  0&  0 & 0 &-1\\
-1& 0&-1& 0 &  0&  0 & 1 & 0
\end{tabular}

$----------------------------------$

V(87):  $\deg = 48$,$\,$ $h^{1,2} = 0$, $rk =1$,
$sq =1$,
$dp =1$, $py =1$,  ${\rm vert}(\Delta^*) =7$:

\begin{tabular}{|c|c|c|c|c|c|c|}
-4& 0&-2&-2& 0 &  0& 1 \\
-1& 0& 0&-1& 1 &  0& 0\\
-1& 0&-1& 0& 0 &  1& 0\\
-1& 1& 0& 0& 0 &  0& 0
\end{tabular}

$----------------------------------$

V(88):  $\deg = 54$,$\,$ $h^{1,2} = 0$, $rk =0$, $sq =0$,
$dp =0$, $py =1$,  ${\rm vert}(\Delta^*) =6$:

\begin{tabular}{|c|c|c|c|c|c|}
-3& 0& 0& 1& 0 & -2 \\
 0& 1& 0& 0& 0 & -1\\
-1& 0& 0& 0& 1 &  0 \\
-1& 0& 1& 0& 0 &  0
\end{tabular}

$----------------------------------$

\newpage 

V(89):  $\deg = 54$,$\,$ $h^{1,2} = 0$, $rk =0$, $sq =0$,
$dp =0$, $py =1$,  ${\rm vert}(\Delta^*) =6$:

\begin{tabular}{|c|c|c|c|c|c|}
-1& 0& 1& 0& 0 & -3 \\
 1& 1& 0& 0& 0 & -1\\
 0& 0& 0& 1& 0 & -1 \\
-1& 0& 0& 0& 1 &  0
\end{tabular}

$----------------------------------$

V(90):  $\deg = 54$,$\,$ $h^{1,2} = 0$, $rk =1$,
$sq =1$,
$dp =1$, $py =1$,  ${\rm vert}(\Delta^*) =7$:

\begin{tabular}{|c|c|c|c|c|c|c|}
-5&-3& 0&-2& 0 &  0& 1 \\
-2&-1& 0&-1& 1 &  0& 0\\
-1&-1& 0& 0& 0 &  1& 0\\
-1& 0& 1& 0& 0 &  0& 0
\end{tabular}

$----------------------------------$

V(91):  $\deg = 56$,$\,$ $h^{1,2} = 0$, $rk =0$, $sq =0$,
$dp =0$, $py =1$,  ${\rm vert}(\Delta^*) =6$:

\begin{tabular}{|c|c|c|c|c|c|}
 2& 0& 0& 1& 0 & -2 \\
 1& 1& 0& 0& 0 & -1\\
-1& 0& 0& 0& 1 &  0 \\
-1& 0& 1& 0& 0 &  0
\end{tabular}

$----------------------------------$

V(92):  $\deg = 62$,$\,$ $h^{1,2} = 0$, $rk =0$, $sq =0$,
$dp =0$, $py =1$,  ${\rm vert}(\Delta^*) =6$:

\begin{tabular}{|c|c|c|c|c|c|}
-1& 0& 0& 1& 0 & -2 \\
 2& 1& 0& 0& 0 & -1\\
-1& 0& 0& 0& 1 &  0 \\
-1& 0& 1& 0& 0 &  0
\end{tabular}

$----------------------------------$

}

\section{$B_2 =3$}

{\tiny

V(93):  $\deg = 12$,  $h^{1,2} = 8$, $rk =5$,
$sq =12$,
$dp =12$, $py =0$,  ${\rm vert}(\Delta^*) =12$:

\begin{tabular}{|c|c|c|c|c|c|c|c|c|c|c|c|}
-1&-1&-1& 1 & -1& 1 & 1 & 1 &-1 & 1 & -1& 1\\
 0& 0& 0&-1 &  0&-1 & 0 & 0 & 1 & 0 &  1& 0\\
 0& 0& 1& 0 &  1& 0 & 0 & 0 & 0 &-1 &  0&-1\\
 0& 1& 1& 0 &  0&-1 &-1 & 0 & 1 & 0 &  0&-1
\end{tabular}

$----------------------------------$

V(94):  $\deg = 18$,  $h^{1,2} = 3$, $rk =3$,
$sq =5$,
$dp =5$, $py =0$,  ${\rm vert}(\Delta^*) =10$:

\begin{tabular}{|c|c|c|c|c|c|c|c|c|c|}
 0& 0& 0& 1 &  0&-1 &-1 & 0 & 0 &-1 \\
 0& 0& 1& 0 &  0& 0 &-1 &-1 & 0 &-1 \\
 0& 1& 0& 0 &  0& 0 &-1 & 0 &-1 & 0 \\
-1& 0& 0& 0 &  1& 0 & 0 & 0 & 0 &-1
\end{tabular}

$----------------------------------$

V(95):  $\deg = 22$,  $h^{1,2} = 1$, $rk =2$,
$sq =2$,
$dp =2$, $py =0$,  ${\rm vert}(\Delta^*) =9$:

\begin{tabular}{|c|c|c|c|c|c|c|c|c|}
 0& 0& 0& 0 & -1&-1&-1 & 0 & 1\\
 0& 0&-1& 0 & -1& 0&-1 & 1 & 0\\
-1& 1& 0& 0 & -1& 0& 0 & 0 & 0\\
 0& 0& 0& 1 &  0& 0&-1 & 0 & 0
\end{tabular}

$----------------------------------$

V(96):  $\deg = 24$,  $h^{1,2} = 1$, $rk =1$,
$sq =1$,
$dp =1$, $py =0$,  ${\rm vert}(\Delta^*) =8$:

\begin{tabular}{|c|c|c|c|c|c|c|c|}
 0& 0& 0&-1 & -1& 0 & 0 & 1\\
 0& 0& 0& 0 & -1&-1 & 1 & 0\\
-1& 0& 1& 0 &  0& 0 & 0 & 0\\
 0& 1& 0& 0 & -1& 0 & 0 & 0
\end{tabular}

$----------------------------------$

V(97):  $\deg = 24$,  $h^{1,2} = 1$, $rk =2$,
$sq =2$,
$dp =2$, $py =0$,  ${\rm vert}(\Delta^*) =9$:

\begin{tabular}{|c|c|c|c|c|c|c|c|c|}
-2& 0&-1& 0 &  0&-1& 0 & 0 & 1\\
-1& 0& 0&-1 &  0&-1& 1 & 0 & 0\\
-1& 0& 0& 0 & -1&-1& 0 & 1 & 0\\
-1& 1& 0& 0 &  0& 0& 0 & 0 & 0
\end{tabular}

$----------------------------------$
\newpage 

V(98):  $\deg = 24$,  $h^{1,2} = 1$, $rk =6$,
$sq =7$,
$dp =7$, $py =1$,  ${\rm vert}(\Delta^*) =13$:

\begin{tabular}{|c|c|c|c|c|c|c|c|c|c|c|c|c|}
-1& 0& 0& 1 & -2& 0 &-1 &-1 & 0 &-3 &-2 &-2& -2 \\
 1& 1& 1& 0 & -1& 0 &-1 & 0 & 0 &-1 & 0 &-1&  0\\
-1&-1& 0& 0 &  1& 0 & 1 &-1 & 1 & 0 & 0 & 0& -1\\
 0& 1& 0& 0 & -1& 1 & 0 & 1 & 0 &-1 &-1 & 0&  0
\end{tabular}

$----------------------------------$

V(99):  $\deg = 26$,  $h^{1,2} = 0$, $rk =5$,
$sq =5$,
$dp =5$, $py =1$,  ${\rm vert}(\Delta^*) =12$:

\begin{tabular}{|c|c|c|c|c|c|c|c|c|c|c|c|}
 1&-1& 0&-1 & -2&-2 & 0 & 1 & 0 &-2 & 0 &-2 \\
 0& 1& 0&-1 & -1& 0 &-1 & 0 & 1 & 0 & 0 & 1\\
 0& 0& 0& 0 &  0&-1 & 1 & 1 & 0 & 0 & 1 &-1\\
 0&-1& 1& 1 &  0& 0 & 1 & 1 & 0 &-1 & 0 &-1
\end{tabular}

$----------------------------------$

V(100):  $\deg = 26$,  $h^{1,2} = 3$, $rk =1$,
$sq =1$,
$dp =3$, $py =0$,  ${\rm vert}(\Delta^*) =8$:

\begin{tabular}{|c|c|c|c|c|c|c|c|}
 0&-1&-2& 0 & -1& 0 & 0 & 1\\
 0& 0&-1&-1 & -2& 1 & 0 & 0\\
 0& 0&-1& 0 & -1& 0 & 1 & 0\\
 1& 0&-1& 0 & -1& 0 & 0 & 0
\end{tabular}

$----------------------------------$

V(101):  $\deg = 28$,  $h^{1,2} = 0$, $rk =4$,
$sq =4$,
$dp =4$, $py =1$,  ${\rm vert}(\Delta^*) =11$:

\begin{tabular}{|c|c|c|c|c|c|c|c|c|c|c|}
 1&-1&-1&-2 &  0& 0 & 1 &-2 & 0 &-2 &-2 \\
 0& 0& 0&-1 &  0& 0 & 1 & 0 & 1 &-1 & 0 \\
 0& 1&-1& 0 &  0& 1 & 0 &-1 & 0 & 1 & 0 \\
 0&-1& 1& 0 &  1& 0 & 1 & 0 & 0 &-1 &-1
\end{tabular}

$----------------------------------$

V(102):  $\deg = 28$,  $h^{1,2} = 0$, $rk =5$,
$sq =5$,
$dp =5$, $py =1$,  ${\rm vert}(\Delta^*) =12$:

\begin{tabular}{|c|c|c|c|c|c|c|c|c|c|c|c|}
-2& 0& 0&-2 & -1&-1 &-2 & 0 &-1 &-2 &-1 & 1 \\
-1& 1& 0& 0 &  0&-1 &-1 & 0 & 0 & 0 &-1 & 0\\
 1& 0& 1&-1 & -1& 0 & 0 & 0 & 1 & 0 & 1 & 0\\
-1& 0& 0& 0 &  1& 1 & 0 & 1 &-1 &-1 & 0 & 0
\end{tabular}

$----------------------------------$

V(103):  $\deg = 28$,  $h^{1,2} = 1$, $rk =0$,
$sq =0$,
$dp =0$, $py =0$,  ${\rm vert}(\Delta^*) =7$:

\begin{tabular}{|c|c|c|c|c|c|c|}
 1& 0& 0& 1&-1 &  0& 0 \\
 0& 1& 0& 0& 0 & -1& 0\\
-1& 0& 0& 0& 0 &  0& 1\\
-1& 0& 1& 0& 0 &  0& 0
\end{tabular}

$----------------------------------$

V(104):  $\deg = 28$,  $h^{1,2} = 1$, $rk =1$,
$sq =1$,
$dp =1$, $py =0$,  ${\rm vert}(\Delta^*) =8$:

\begin{tabular}{|c|c|c|c|c|c|c|c|}
-2& 0& 0&-1 & -1& 0 & 0 & 1\\
-1& 0&-1&-1 &  0& 1 & 0 & 0\\
-1& 0& 0&-1 &  0& 0 & 1 & 0\\
-1& 1& 0& 0 &  0& 0 & 0 & 0
\end{tabular}

$----------------------------------$

V(105):  $\deg = 28$,  $h^{1,2} = 1$, $rk =4$,
$sq =5$,
$dp =5$, $py =1$,  ${\rm vert}(\Delta^*) =11$:

\begin{tabular}{|c|c|c|c|c|c|c|c|c|c|c|}
 1& 1& 0& 0 &  0& 0 &-2 &-2 &-1 &-1 &-2 \\
 0& 1& 1& 0 &  0& 1 & 0 &-1 &-1 & 0 & 0 \\
 0& 1& 0& 0 &  1& 1 &-1 & 0 & 0 &-1 & 0 \\
 0& 0& 0& 1 &  0&-1 & 0 & 0 & 1 & 1 &-1
\end{tabular}

$----------------------------------$

V(106):  $\deg = 28$,  $h^{1,2} = 1$, $rk =4$,
$sq =5$,
$dp =5$, $py =1$,  ${\rm vert}(\Delta^*) =11$:

\begin{tabular}{|c|c|c|c|c|c|c|c|c|c|c|}
-2&-1&-1&-3 &  0&-1 &-2 &-2 & 0 & 1 & 0 \\
-1&-1&-1& 0 &  1& 1 & 1 & 0 & 0 & 0 & 0 \\
 0& 0& 1&-1 &  0&-1 &-1 &-1 & 0 & 0 & 1 \\
 0& 1& 0&-1 &  0& 0 &-1 & 0 & 1 & 0 & 0
\end{tabular}

$----------------------------------$

V(107):  $\deg = 30$,  $h^{1,2} = 0$, $rk =3$,
$sq =3$,
$dp =3$, $py =1$,  ${\rm vert}(\Delta^*) =10$:

\begin{tabular}{|c|c|c|c|c|c|c|c|c|c|}
 1& 1& 0& 0 &  0&-3 &-2 &-1 &-2 &-2 \\
 0& 1& 1& 0 &  0& 0 &-1 &-1 & 0 & 0 \\
 0& 1& 0& 0 &  1&-1 & 0 & 0 &-1 & 0 \\
 0& 0& 0& 1 &  0&-1 & 0 & 1 & 0 &-1
\end{tabular}

$----------------------------------$

\newpage 

V(108):  $\deg = 30$,  $h^{1,2} = 0$, $rk =4$,
$sq =4$,
$dp =4$, $py =1$,  ${\rm vert}(\Delta^*) =11$:

\begin{tabular}{|c|c|c|c|c|c|c|c|c|c|c|}
-1& 0& 0& 1 &  0&-3 &-2 &-2 &-1 &-2 & 0 \\
-1&-1& 0& 0 &  1& 0 & 0 & 0 &-1 &-1 & 0 \\
 1& 1& 1& 0 &  0&-1 & 0 &-1 & 0 & 0 & 0 \\
 0& 1& 0& 0 &  0&-1 &-1 & 0 & 1 & 0 & 1
\end{tabular}

$----------------------------------$

V(109):  $\deg = 30$,  $h^{1,2} = 0$, $rk =4$,
$sq =4$,
$dp =4$, $py =1$,  ${\rm vert}(\Delta^*) =11$:

\begin{tabular}{|c|c|c|c|c|c|c|c|c|c|c|}
-3&-2&-1&-3 & -2&-4 & 0 & 0 &-1 & 0 & 1 \\
-1&-1& 1& 0 &  0&-1 & 1 & 0 &-1 & 0 & 0 \\
 0& 0&-1&-1 & -1&-1 & 0 & 0 & 1 & 1 & 0 \\
-1& 0& 0&-1 &  0&-1 & 0 & 1 & 0 & 0 & 0
\end{tabular}

$----------------------------------$

V(110):  $\deg = 32$,  $h^{1,2} = 0$, $rk =3$,
$sq =3$,
$dp =3$, $py =1$,  ${\rm vert}(\Delta^*) =10$:

\begin{tabular}{|c|c|c|c|c|c|c|c|c|c|}
-1&-2&-2& 0 &  0& 1 & 0 &-2 &-1 &-2 \\
 0&-1& 0& 0 &  1& 0 & 0 & 0 & 1 &-1 \\
-1&-1&-1& 0 &  0& 0 & 1 & 0 & 0 & 0 \\
 1& 1& 0& 1 &  0& 0 & 0 &-1 &-1 & 0
\end{tabular}

$----------------------------------$

V(111):  $\deg = 32$,  $h^{1,2} = 0$, $rk =3$,
$sq =3$,
$dp =3$, $py =1$,  ${\rm vert}(\Delta^*) =10$:

\begin{tabular}{|c|c|c|c|c|c|c|c|c|c|}
-2& 0&-1&-2 &  0&-3 &-2 &-1 & 0 & 1 \\
 0& 0&-1&-1 &  0&-1 & 0 & 0 & 1 & 0 \\
-1& 0& 1& 0 &  1& 0 & 0 &-1 & 0 & 0 \\
 0& 1& 0& 0 &  0&-1 &-1 & 1 & 0 & 0
\end{tabular}

$----------------------------------$

V(112):  $\deg = 32$,  $h^{1,2} = 0$, $rk =3$,
$sq =3$,
$dp =3$, $py =1$,  ${\rm vert}(\Delta^*) =10$:

\begin{tabular}{|c|c|c|c|c|c|c|c|c|c|}
-3& 0&-2&-1 & -2& 0 &-4 &-1 & 1 & 0 \\
 0& 1& 0&-1 & -1& 0 &-1 & 1 & 0 & 0 \\
-1& 0&-1& 1 &  0& 0 &-1 &-1 & 0 & 1 \\
-1& 0& 0& 0 &  0& 1 &-1 & 0 & 0 & 0
\end{tabular}

$----------------------------------$

V(113):  $\deg = 32$,  $h^{1,2} = 0$, $rk =4$,
$sq =4$,
$dp =4$, $py =1$,  ${\rm vert}(\Delta^*) =11$:

\begin{tabular}{|c|c|c|c|c|c|c|c|c|c|c|}
-2& 0&-2& 0 & -3&-2 &-3 &-2 &-1 & 0 & 1 \\
 0& 0&-1& 0 & -1&-1 &-1 & 0 & 0 & 1 & 0 \\
-1& 0& 0& 1 &  0&-1 &-1 & 0 &-1 & 0 & 0 \\
 0& 1& 0& 0 & -1& 1 & 0 &-1 & 1 & 0 & 0
\end{tabular}

$----------------------------------$

V(114):  $\deg = 32$,  $h^{1,2} = 1$, $rk =0$,
$sq =0$,
$dp =0$, $py =0$,  ${\rm vert}(\Delta^*) =7$:

\begin{tabular}{|c|c|c|c|c|c|c|}
-1& 0&-1& 0& 0 &  0& 1 \\
-2& 0& 0&-1& 1 &  0& 0\\
-1& 0& 0& 0& 0 &  1& 0\\
-1& 1& 0& 0& 0 &  0& 0
\end{tabular}

$----------------------------------$

V(115):  $\deg = 32$,  $h^{1,2} = 1$, $rk =2$,
$sq =3$,
$dp =3$, $py =1$,  ${\rm vert}(\Delta^*) =9$:

\begin{tabular}{|c|c|c|c|c|c|c|c|c|}
 1&-2&-3&-2 & -1& 0& 0 & 0 & 1\\
 0& 1& 0&-1 & -1& 0& 1 & 0 & 1\\
 0&-1&-1& 0 &  0& 1& 0 & 0 & 1\\
 0&-1&-1& 0 &  1& 0& 0 & 1 & 0
\end{tabular}

$----------------------------------$

V(116):  $\deg = 34$,  $h^{1,2} = 0$, $rk =2$,
$sq =2$,
$dp =2$, $py =1$,  ${\rm vert}(\Delta^*) =9$:

\begin{tabular}{|c|c|c|c|c|c|c|c|c|}
-1& 0& 0& 1 &  0&-2&-1 &-2 &-2\\
 0& 0& 1& 0 &  0& 0& 1 &-1 &-1\\
 1& 1& 0& 0 &  0&-1&-1 & 0 & 1\\
-1& 0& 0& 0 &  1& 0& 0 & 0 &-1
\end{tabular}

$----------------------------------$

V(117):  $\deg = 34$,  $h^{1,2} = 0$, $rk =3$,
$sq =3$,
$dp =3$, $py =1$,  ${\rm vert}(\Delta^*) =10$:

\begin{tabular}{|c|c|c|c|c|c|c|c|c|c|}
 0&-1&-2&-2 & -1& 0 &-2 &-3 & 0 & 1 \\
 0&-1&-1& 0 &  0& 0 &-1 &-1 & 1 & 0 \\
 0& 1& 0&-1 & -1& 1 & 1 & 0 & 0 & 0 \\
 1& 0& 0& 0 &  1& 0 &-1 &-1 & 0 & 0
\end{tabular}

$----------------------------------$
\newpage 

V(118):  $\deg = 34$,  $h^{1,2} = 0$, $rk =3$,
$sq =3$,
$dp =3$, $py =1$,  ${\rm vert}(\Delta^*) =10$:

\begin{tabular}{|c|c|c|c|c|c|c|c|c|c|}
-2&-1&-3& 0 & -2& 0 & 0 &-3 &-2 & 1 \\
-1& 0&-1& 0 &  0& 1 & 0 &-1 &-1 & 0 \\
 1& 1& 0& 1 & -1& 0 & 0 &-1 & 0 & 0 \\
-1&-1&-1& 0 &  0& 0 & 1 & 0 & 0 & 0
\end{tabular}

$----------------------------------$

V(119):  $\deg = 36$,  $h^{1,2} = 0$, $rk =2$,
$sq =2$,
$dp =2$, $py =1$,  ${\rm vert}(\Delta^*) =9$:

\begin{tabular}{|c|c|c|c|c|c|c|c|c|}
-1& 0& 0& 1 &  0&-3&-2 &-2 &-2\\
 0& 0& 1& 0 &  0&-1&-1 &-1 & 0\\
 1& 1& 0& 0 &  0&-1& 0 & 1 &-1\\
-1& 0& 0& 0 &  1& 0& 0 &-1 & 0
\end{tabular}

$----------------------------------$

V(120):  $\deg = 36$,  $h^{1,2} = 0$, $rk =2$,
$sq =2$,
$dp =2$, $py =1$,  ${\rm vert}(\Delta^*) =9$:

\begin{tabular}{|c|c|c|c|c|c|c|c|c|}
-2& 0&-2&-1 &  0&-2 &-3 & 0& 1\\
-1& 0& 0& 0 &  0&-1 &-1 & 1& 0\\
 0& 1&-1&-1 &  0& 1 & 0 & 0& 0\\
 0& 0& 0& 1 &  1&-1 &-1 & 0& 0
\end{tabular}

$----------------------------------$

V(121):  $\deg = 36$,  $h^{1,2} = 0$, $rk =2$,
$sq =2$,
$dp =2$, $py =1$,  ${\rm vert}(\Delta^*) =9$:

\begin{tabular}{|c|c|c|c|c|c|c|c|c|}
 1& 0&-2&-2 &  0& 1 & 0 &-3&-1\\
 0& 0&-1& 0 &  0& 1 & 1 & 0& 1\\
 0& 0& 0&-1 &  1& 1 & 0 &-1& 0\\
 0& 1& 0& 0 &  0& 0 & 0 &-1&-1
\end{tabular}

$----------------------------------$

V(122):  $\deg = 36$,  $h^{1,2} = 0$, $rk =3$,
$sq =3$,
$dp =3$, $py =1$,  ${\rm vert}(\Delta^*) =10$:

\begin{tabular}{|c|c|c|c|c|c|c|c|c|c|}
 1& 0&-2& 0 &  0& 0 & 1 &-2 &-1 &-1 \\
 0& 1&-1& 0 &  1& 0 & 1 & 0 & 0 & 1 \\
 0& 1& 0& 1 &  0& 0 & 1 &-1 & 1 & 0 \\
 0&-1& 0& 0 &  0& 1 & 0 & 0 &-1 &-1
\end{tabular}

$----------------------------------$

V(123):  $\deg = 36$,  $h^{1,2} = 0$, $rk =3$,
$sq =3$,
$dp =3$, $py =1$,  ${\rm vert}(\Delta^*) =10$:

\begin{tabular}{|c|c|c|c|c|c|c|c|c|c|}
-3& 0&-1& 0 & -3&-2 &-4 &-2 & 0 & 1 \\
-1& 0& 0& 0 & -1& 0 &-1 &-1 & 1 & 0 \\
 0& 1&-1& 0 & -1&-1 &-1 & 0 & 0 & 0 \\
-1& 0& 1& 1 &  0& 0 &-1 & 0 & 0 & 0
\end{tabular}

$----------------------------------$

V(124):  $\deg = 38$,  $h^{1,2} = 0$, $rk =1$,
$sq =1$,
$dp =1$, $py =1$,  ${\rm vert}(\Delta^*) =8$:

\begin{tabular}{|c|c|c|c|c|c|c|c|}
-1&-2&-3& 0 &  1&  0 &-2 & 0\\
-1&-1& 0& 0 &  0&  0 & 1 & 1\\
 0& 0&-1& 0 &  0&  1 &-1 & 0\\
 1& 0&-1& 1 &  0&  0 &-1 & 0
\end{tabular}

$----------------------------------$

V(125):  $\deg = 38$,  $h^{1,2} = 0$, $rk =1$,
$sq =1$,
$dp =1$, $py =1$,  ${\rm vert}(\Delta^*) =8$:

\begin{tabular}{|c|c|c|c|c|c|c|c|}
 1&-2&-2& 0 &  0&  0 & 1 &-4\\
 0&-1& 0& 1 &  0&  0 & 1 &-1\\
 0& 0&-1& 0 &  0&  1 & 1 &-1\\
 0& 0& 0& 0 &  1&  0 & 0 &-1
\end{tabular}

$----------------------------------$

V(126):  $\deg = 38$,  $h^{1,2} = 0$, $rk =2$,
$sq =2$,
$dp =2$, $py =1$,  ${\rm vert}(\Delta^*) =9$:

\begin{tabular}{|c|c|c|c|c|c|c|c|c|}
 1& 1& 0& 0 &  0&-2& 0 &-2 &-2\\
 0& 0& 0& 1 &  0& 0&-1 &-1 & 0\\
 0& 1& 0& 0 &  1&-1& 1 & 0 & 0\\
 0& 1& 1& 0 &  0& 0& 1 & 0 &-1
\end{tabular}

$----------------------------------$

V(127):  $\deg = 38$,  $h^{1,2} = 0$, $rk =2$,
$sq =2$,
$dp =2$, $py =1$,  ${\rm vert}(\Delta^*) =9$:

\begin{tabular}{|c|c|c|c|c|c|c|c|c|}
-1& 0& 0& 1 &  0& 0&-2 &-2 &-1\\
 1& 1& 1& 0 &  0& 0&-1 & 0 & 0\\
 0& 1& 0& 0 &  0& 1& 0 &-1 & 1\\
-1&-1& 0& 0 &  1& 0& 0 & 0 &-1
\end{tabular}

$----------------------------------$

\newpage 

V(128):  $\deg = 38$,  $h^{1,2} = 0$, $rk =2$,
$sq =2$,
$dp =2$, $py =1$,  ${\rm vert}(\Delta^*) =9$:

\begin{tabular}{|c|c|c|c|c|c|c|c|c|}
 1&-2&-1&-3 &  0&-2& 1 & 0 & 0\\
 0& 1& 1& 0 &  0&-1& 1 & 1 & 0\\
 0&-1&-1&-1 &  1& 0& 0 & 0 & 0\\
 0&-1& 0&-1 &  0& 0& 1 & 0 & 1
\end{tabular}

$----------------------------------$

V(129):  $\deg = 38$,  $h^{1,2} = 0$, $rk =2$,
$sq =2$,
$dp =2$, $py =1$,  ${\rm vert}(\Delta^*) =9$:

\begin{tabular}{|c|c|c|c|c|c|c|c|c|}
-2&-1& 0&-2 &  0&-3 &-1 & 0& 1\\
 0&-1& 0&-1 &  0&-1 & 0 & 1& 0\\
-1& 1& 0& 0 &  1& 0 &-1 & 0& 0\\
 0& 0& 1& 0 &  0&-1 & 1 & 0& 0
\end{tabular}

$----------------------------------$

V(130):  $\deg = 38$,  $h^{1,2} = 0$, $rk =3$,
$sq =3$,
$dp =3$, $py =1$,  ${\rm vert}(\Delta^*) =10$:

\begin{tabular}{|c|c|c|c|c|c|c|c|c|c|}
 1& 1& 0&-2 & -1& 0 & 0 & 0 &-2 &-3 \\
 0& 0& 0&-1 & -1&-1 & 0 & 1 & 0 &-1 \\
 0& 1& 0& 0 &  1& 1 & 1 & 0 &-1 & 0 \\
 0& 1& 1& 0 &  0& 1 & 0 & 0 & 0 &-1
\end{tabular}

$----------------------------------$

V(131):  $\deg = 40$,  $h^{1,2} = 0$, $rk =1$,
$sq =1$,
$dp =1$, $py =1$,  ${\rm vert}(\Delta^*) =8$:

\begin{tabular}{|c|c|c|c|c|c|c|c|}
 1&-2&-2& 0 &  1&  0 &-3 & 0\\
 0& 1&-1& 0 &  1&  1 & 0 & 0\\
 0&-1& 0& 0 &  0&  0 &-1 & 1\\
 0&-1& 0& 1 &  1&  0 &-1 & 0
\end{tabular}

$----------------------------------$

V(132):  $\deg = 42$,  $h^{1,2} = 0$, $rk =1$,
$sq =1$,
$dp =1$, $py =1$,  ${\rm vert}(\Delta^*) =8$:

\begin{tabular}{|c|c|c|c|c|c|c|c|}
 1& 1& 0&-1 &  0& -2 & 0 &-3\\
 0& 1& 1& 1 &  0& -1 & 0 & 0\\
 0& 0& 0&-1 &  0&  0 & 1 &-1\\
 0& 1& 0& 0 &  1&  0 & 0 &-1
\end{tabular}

$----------------------------------$

V(133):  $\deg = 42$,  $h^{1,2} = 0$, $rk =1$,
$sq =1$,
$dp =1$, $py =1$,  ${\rm vert}(\Delta^*) =8$:

\begin{tabular}{|c|c|c|c|c|c|c|c|}
-2& 0&-2& 0 & -3& -1 & 1 & 0\\
-1& 0& 0& 1 &  0& -1 & 0 & 0\\
 0& 0&-1& 0 & -1&  0 & 0 & 1\\
 0& 1& 0& 0 & -1&  1 & 0 & 0
\end{tabular}

$----------------------------------$

V(134):  $\deg = 42$,  $h^{1,2} = 0$, $rk =1$,
$sq =1$,
$dp =1$, $py =1$,  ${\rm vert}(\Delta^*) =8$:

\begin{tabular}{|c|c|c|c|c|c|c|c|}
-4&-2&-2 &  0&  0 &-1 & 0 & 1\\
-1&-1& 0 &  1&  0 &-1 & 0 & 0\\
-1& 0&-1 &  0&  0 & 1 & 1 & 0\\
-1& 0& 0 &  0&  1 & 0 & 0 & 0
\end{tabular}

$----------------------------------$

V(135):  $\deg = 42$,  $h^{1,2} = 0$, $rk =2$,
$sq =2$,
$dp =2$, $py =1$,  ${\rm vert}(\Delta^*) =9$:

\begin{tabular}{|c|c|c|c|c|c|c|c|c|}
 1& 2& 1&-2 &  0& 0 & 0 &-3&-3\\
 0& 1& 0&-1 &  0& 0 & 1 &-1&-1\\
 0& 1& 1& 0 &  0& 1 & 0 &-1& 0\\
 0& 1& 1& 0 &  1& 0 & 0 & 0&-1
\end{tabular}

$----------------------------------$

V(136):  $\deg = 42$,  $h^{1,2} = 0$, $rk =2$,
$sq =2$,
$dp =2$, $py =1$,  ${\rm vert}(\Delta^*) =9$:

\begin{tabular}{|c|c|c|c|c|c|c|c|c|}
-4&-1& 0&-2 & -3&-2 & 0 & 0& 1\\
-1& 0& 0&-1 & -1& 0 & 1 & 0& 0\\
-1&-1& 0& 0 & -1&-1 & 0 & 1& 0\\
-1& 1& 1& 0 &  0& 0 & 0 & 0& 0
\end{tabular}

$----------------------------------$

V(137):  $\deg = 44$,  $h^{1,2} = 0$, $rk =0$,
$sq =0$,
$dp =0$, $py =1$,  ${\rm vert}(\Delta^*) =7$:

\begin{tabular}{|c|c|c|c|c|c|c|}
-1& 0& 0& 1& 0 & -3&-2 \\
 1& 1& 0& 0& 0 & -1&-1\\
 0& 0& 1& 0& 0 & -1& 0\\
-1& 0& 0& 0& 1 &  0& 0
\end{tabular}

$----------------------------------$

\newpage 

V(138):  $\deg = 44$,  $h^{1,2} = 0$, $rk =1$,
$sq =1$,
$dp =1$, $py =1$,  ${\rm vert}(\Delta^*) =8$:

\begin{tabular}{|c|c|c|c|c|c|c|c|}
-3& 0& 0 & -3& -4 &-2 & 0 & 1\\
-1& 0& 0 & -1& -1 &-1 & 1 & 0\\
 0& 1& 0 & -1& -1 & 0 & 0 & 0\\
-1& 0& 1 &  0& -1 & 0 & 0 & 0
\end{tabular}

$----------------------------------$

V(139):  $\deg = 46$,  $h^{1,2} = 0$, $rk =0$,
$sq =0$,
$dp =0$, $py =1$,  ${\rm vert}(\Delta^*) =7$:

\begin{tabular}{|c|c|c|c|c|c|c|}
 1& 0&-2& 0& 1 &  0&-3 \\
 0& 1&-1& 0& 1 &  0& 0\\
 0& 0& 0& 0& 1 &  1&-1\\
 0& 0& 0& 1& 0 &  0&-1
\end{tabular}

$----------------------------------$

V(140):  $\deg = 46$,  $h^{1,2} = 0$, $rk =1$,
$sq =1$,
$dp =1$, $py =1$,  ${\rm vert}(\Delta^*) =8$:

\begin{tabular}{|c|c|c|c|c|c|c|c|}
-4& 0&-3 & -1& -2 & 0 & 0 & 1\\
-1& 0&-1 & -1& -1 & 1 & 0 & 0\\
-1& 0&-1 &  0&  0 & 0 & 1 & 0\\
-1& 1& 0 &  1&  0 & 0 & 0 & 0
\end{tabular}

$----------------------------------$

V(141):  $\deg = 48$,  $h^{1,2} = 0$, $rk =0$,
$sq =0$,
$dp =0$, $py =1$,  ${\rm vert}(\Delta^*) =7$:

\begin{tabular}{|c|c|c|c|c|c|c|}
-2& 0& 0& 1& 0 & -2&-2 \\
 0& 0& 1& 0& 0 & -1& 0\\
 0& 1& 0& 0& 0 &  0&-1\\
-1& 0& 0& 0& 1 &  0& 0
\end{tabular}

$----------------------------------$

V(142):  $\deg = 48$,  $h^{1,2} = 0$, $rk =0$,
$sq =0$,
$dp =0$, $py =1$,  ${\rm vert}(\Delta^*) =7$:

\begin{tabular}{|c|c|c|c|c|c|c|}
-1& 0& 0& 1& 0 & -2&-2 \\
 1& 0& 1& 0& 0 & -1& 0\\
 0& 1& 0& 0& 0 &  0&-1\\
-1& 0& 0& 0& 1 &  0& 0
\end{tabular}

$----------------------------------$

V(143):  $\deg = 48$,  $h^{1,2} = 0$, $rk =1$,
$sq =1$,
$dp =1$, $py =1$,  ${\rm vert}(\Delta^*) =8$:

\begin{tabular}{|c|c|c|c|c|c|c|c|}
-3& 0& 0&-2 & -2& -1 & 0 & 1\\
-2& 0& 0&-1 &  0& -1 & 1 & 0\\
 1& 1& 0& 0 & -1&  1 & 0 & 0\\
-1& 0& 1& 0 &  0&  0 & 0 & 0
\end{tabular}

$----------------------------------$

V(144):  $\deg = 50$,  $h^{1,2} = 0$, $rk =0$,
$sq =0$,
$dp =0$, $py =1$,  ${\rm vert}(\Delta^*) =7$:

\begin{tabular}{|c|c|c|c|c|c|c|}
-1& 0& 0& 1& 0 & -2&-1 \\
 0& 0& 1& 0& 0 & -1& 1\\
 1& 1& 0& 0& 0 &  0&-1\\
-1& 0& 0& 0& 1 &  0& 0
\end{tabular}

$----------------------------------$

V(145):  $\deg = 50$,  $h^{1,2} = 0$, $rk =0$,
$sq =0$,
$dp =0$, $py =1$,  ${\rm vert}(\Delta^*) =7$:

\begin{tabular}{|c|c|c|c|c|c|c|}
 1& 0&-2& 0& 1 &  0&-2 \\
 0& 1&-1& 0& 1 &  0& 1\\
 0& 0& 0& 0& 1 &  1&-1\\
 0& 0& 0& 1& 0 &  0&-1
\end{tabular}

$----------------------------------$

V(146):  $\deg = 52$,  $h^{1,2} = 0$, $rk =0$,
$sq =0$,
$dp =0$, $py =1$,  ${\rm vert}(\Delta^*) =7$:

\begin{tabular}{|c|c|c|c|c|c|c|}
-1& 0& 0& 1& 0 & -2&-1 \\
 1& 0& 1& 0& 0 & -1& 1\\
 0& 1& 0& 0& 0 &  0&-1\\
-1& 0& 0& 0& 1 &  0& 0
\end{tabular}

$----------------------------------$
}

\newpage 
\section{$B_2 =4$}

{\tiny

V(147):  $\deg = 24$,   $h^{1,2} = 1$, $rk =0$,
$sq =0$,
$dp =0$, $py =0$, ${\rm vert}(\Delta^*) =8$:

\begin{tabular}{|c|c|c|c|c|c|c|c|}
 0& 0& 0& 1 &  0&-1 & 0 & 0\\
 0& 0& 1& 0 &  0& 0 &-1 & 0\\
 0& 1& 0& 0 &  0& 0 & 0 &-1\\
-1& 0& 0& 0 &  1& 0 & 0 & 0
\end{tabular}

$----------------------------------$

V(148):  $\deg = 24$,   $h^{1,2} = 1$, $rk =5$,
$sq =6$,
$dp =6$, $py =1$,  ${\rm vert}(\Delta^*) =13$:

\begin{tabular}{|c|c|c|c|c|c|c|c|c|c|c|c|c|}
-1&-2&-2& 0 &  0& 1 &-1 & 0 &-2 &-1 &-2 & 0& -1 \\
 1& 0& 0& 1 &  1& 0 & 0 & 0 &-1 &-1 &-1 & 0&  0\\
 0& 0&-1& 1 &  0& 0 & 1 & 1 & 0 & 0 &-1 & 0& -1\\
-1&-1& 0&-1 &  0& 0 &-1 & 0 & 0 & 1 & 1 & 1&  1
\end{tabular}

$----------------------------------$

V(149):  $\deg = 28$,   $h^{1,2} = 1$, $rk =0$,
$sq =0$,
$dp =0$, $py =0$,  ${\rm vert}(\Delta^*) =8$:

\begin{tabular}{|c|c|c|c|c|c|c|c|}
-1& 0& 0&-1 & -1& 0 & 0 & 1\\
-1& 0& 0& 0 & -1&-1 & 1 & 0\\
 0& 1& 0& 0 & -1& 0 & 0 & 0\\
-1& 0& 1& 0 &  0& 0 & 0 & 0
\end{tabular}

$----------------------------------$

V(150):  $\deg = 28$,   $h^{1,2} = 1$, $rk =3$,
$sq =4$,
$dp =4$, $py =1$,  ${\rm vert}(\Delta^*) =11$:

\begin{tabular}{|c|c|c|c|c|c|c|c|c|c|c|}
-1&-2&-2& 0 &  0& 1 & 0 &-2 &-1 & 0 &-2 \\
 1& 0& 0& 1 &  1& 0 & 0 &-1 &-1 & 0 &-1 \\
 0& 0&-1& 1 &  0& 0 & 0 &-1 & 0 & 1 & 0 \\
-1&-1& 0&-1 &  0& 0 & 1 & 1 & 1 & 0 & 0
\end{tabular}

$----------------------------------$

V(151):  $\deg = 30$,   $h^{1,2} = 0$, $rk =3$,
$sq =3$,
$dp =3$, $py =1$,  ${\rm vert}(\Delta^*) =11$:

\begin{tabular}{|c|c|c|c|c|c|c|c|c|c|c|}
-2& 0&-1&-1 &  0& 1 & 0 &-1 &-1 &-2 &-2 \\
-1& 0&-1&-1 &  0& 0 & 1 & 0 & 0 & 0 & 0 \\
 0& 0& 0& 1 &  1& 0 & 0 & 1 &-1 &-1 & 0 \\
 0& 1& 1& 0 &  0& 0 & 0 &-1 & 1 & 0 &-1
\end{tabular}

$----------------------------------$

V(152):  $\deg = 32$,   $h^{1,2} = 0$, $rk =2$,
$sq =2$,
$dp =2$, $py =1$,  ${\rm vert}(\Delta^*) =10$:

\begin{tabular}{|c|c|c|c|c|c|c|c|c|c|}
-2& 0& 0&-1 &  1&-1 & 0 &-1 &-2 &-2 \\
 0& 1& 0& 0 &  0&-1 & 0 & 0 & 0 &-1 \\
 0& 0& 1& 1 &  0& 1 & 0 &-1 &-1 & 0 \\
-1& 0& 0&-1 &  0& 0 & 1 & 1 & 0 & 0
\end{tabular}

$----------------------------------$

V(153):  $\deg = 32$,   $h^{1,2} = 0$, $rk =2$,
$sq =2$,
$dp =2$, $py =1$,  ${\rm vert}(\Delta^*) =10$:

\begin{tabular}{|c|c|c|c|c|c|c|c|c|c|}
-1& 0& 0&-1 & -2& 1 &-2 &-1 & 0 &-1 \\
 1& 1& 0& 1 &  0& 0 &-1 &-1 & 0 & 0 \\
 0& 0& 0&-1 & -1& 0 & 0 & 1 & 1 & 1 \\
-1& 0& 1& 0 &  0& 0 & 0 & 0 & 0 &-1
\end{tabular}

$----------------------------------$

V(154):  $\deg = 34$,   $h^{1,2} = 0$, $rk =2$,
$sq =2$,
$dp =2$, $py =1$,  ${\rm vert}(\Delta^*) =10$:

\begin{tabular}{|c|c|c|c|c|c|c|c|c|c|}
-2& 0& 0&-2 & -3&-2 &-1 &-3 & 0 & 1 \\
 0& 0& 0&-1 & -1& 0 & 0 &-1 & 1 & 0 \\
-1& 0& 1& 0 &  0& 0 &-1 &-1 & 0 & 0 \\
 0& 1& 0& 0 & -1&-1 & 1 & 0 & 0 & 0
\end{tabular}

$----------------------------------$

V(155):  $\deg = 34$,   $h^{1,2} = 0$, $rk =3$,
$sq =3$,
$dp =3$, $py =1$,  ${\rm vert}(\Delta^*) =11$:

\begin{tabular}{|c|c|c|c|c|c|c|c|c|c|c|}
-2&-3& 0& 1 &  0&-2 & 0 &-3 &-3 &-2 &-4 \\
 0& 0& 1& 0 &  0&-1 & 0 &-1 &-1 & 0 &-1 \\
 0&-1& 0& 0 &  1& 0 & 0 &-1 & 0 &-1 &-1 \\
-1&-1& 0& 0 &  0& 0 & 1 & 0 &-1 & 0 &-1
\end{tabular}

$----------------------------------$

V(156):  $\deg = 36$,   $h^{1,2} = 0$, $rk =1$,
$sq =1$,
$dp =1$, $py =1$,  ${\rm vert}(\Delta^*) =9$:

\begin{tabular}{|c|c|c|c|c|c|c|c|c|}
-2& 0& 0&-2 & -3&-2&-1 & 0 & 1\\
-1& 0& 1& 0 &  0& 0&-1 & 0 & 0\\
 0& 0& 0&-1 & -1& 0& 0 & 1 & 0\\
 0& 1& 0& 0 & -1&-1& 1 & 0 & 0
\end{tabular}

$----------------------------------$

V(157):  $\deg = 36$,   $h^{1,2} = 0$, $rk =1$,
$sq =1$,
$dp =1$, $py =1$,  ${\rm vert}(\Delta^*) =9$:

\begin{tabular}{|c|c|c|c|c|c|c|c|c|}
-4&-2&-1&-2 &  0& 0 &-1 & 0& 1\\
-1&-1& 1& 0 &  1& 0 &-1 & 0& 0\\
-1& 0&-1&-1 &  0& 0 & 1 & 1& 0\\
-1& 0& 0& 0 &  0& 1 & 0 & 0& 0
\end{tabular}

$----------------------------------$

V(158):  $\deg = 36$,   $h^{1,2} = 0$, $rk =2$,
$sq =2$,
$dp =2$, $py =1$,  ${\rm vert}(\Delta^*) =10$:

\begin{tabular}{|c|c|c|c|c|c|c|c|c|c|}
-2& 0&-2& 0 & -1&-2 &-1 &-2 & 0 & 1 \\
 0& 0&-1& 0 & -1& 0 & 0 &-1 & 1 & 0 \\
-1& 0& 0& 1 &  0& 0 &-1 &-1 & 0 & 0 \\
 0& 1& 0& 0 &  1&-1 & 1 & 1 & 0 & 0
\end{tabular}

$----------------------------------$

V(159):  $\deg = 38$,   $h^{1,2} = 0$, $rk =1$,
$sq =1$,
$dp =1$, $py =1$,  ${\rm vert}(\Delta^*) =9$:

\begin{tabular}{|c|c|c|c|c|c|c|c|c|}
-2& 0&-1& 0 &  1& 0&-1 &-2 &-2\\
 0& 0& 0& 1 &  0& 0&-1 &-1 & 0\\
-1& 0&-1& 0 &  0& 1& 0 & 0 & 0\\
 0& 1& 1& 0 &  0& 0& 1 & 0 &-1
\end{tabular}

$----------------------------------$

V(160):  $\deg = 40$,   $h^{1,2} = 0$, $rk =0$,
$sq =0$,
$dp =0$, $py =1$,  ${\rm vert}(\Delta^*) =8$:

\begin{tabular}{|c|c|c|c|c|c|c|c|}
 1& 0&-2&-2 &  0&  1 & 0 &-3\\
 0& 1& 0&-1 &  0&  1 & 0 & 0\\
 0& 0&-1& 0 &  0&  1 & 1 &-1\\
 0& 0& 0& 0 &  1&  0 & 0 &-1
\end{tabular}

$----------------------------------$

V(161):  $\deg = 40$,   $h^{1,2} = 0$, $rk =1$,
$sq =1$,
$dp =1$, $py =1$,  ${\rm vert}(\Delta^*) =9$:

\begin{tabular}{|c|c|c|c|c|c|c|c|c|}
 0&-1&-2&-2 &  0&-2&-3 & 0 & 1\\
 0&-1&-1& 0 &  0&-1&-1 & 1 & 0\\
 0& 1& 0&-1 &  1& 1& 0 & 0 & 0\\
 1& 0& 0& 0 &  0&-1&-1 & 0 & 0
\end{tabular}

$----------------------------------$

V(162):  $\deg = 42$,   $h^{1,2} = 0$, $rk =0$,
$sq =0$,
$dp =0$, $py =1$,  ${\rm vert}(\Delta^*) =8$:

\begin{tabular}{|c|c|c|c|c|c|c|c|}
 1& 0& 1& 0 &  0&-2 &-2 &-2\\
 0& 0& 1& 1 &  0&-1 & 0 & 0\\
 0& 0& 1& 0 &  1& 0 &-1 & 0\\
 0& 1& 0& 0 &  0& 0 & 0 &-1
\end{tabular}

$----------------------------------$

V(163):  $\deg = 44$,   $h^{1,2} = 0$, $rk =0$,
$sq =0$,
$dp =0$, $py =1$,  ${\rm vert}(\Delta^*) =8$:

\begin{tabular}{|c|c|c|c|c|c|c|c|}
 1& 0& 0& 0 &  1&-2 &-2 &-1\\
 0& 0& 0& 1 &  1&-1 & 0 & 1\\
 0& 0& 1& 0 &  1& 0 &-1 & 0\\
 0& 1& 0& 0 &  0& 0 & 0 &-1
\end{tabular}

$----------------------------------$

V(164):  $\deg = 46$,   $h^{1,2} = 0$, $rk =0$,
$sq =0$,
$dp =0$, $py =1$,  ${\rm vert}(\Delta^*) =8$:

\begin{tabular}{|c|c|c|c|c|c|c|c|}
 1& 0& 1& 0 &  0&-2 &-2 & 0\\
 0& 0& 1& 1 &  0&-1 & 0 & 1\\
 0& 0& 1& 0 &  1& 0 &-1 & 1\\
 0& 1& 0& 0 &  0& 0 & 0 &-1
\end{tabular}

$----------------------------------$

}
\newpage 
\section{$B_2 =5$}

{\tiny

V(165):  $\deg = 36$,  $h^{1,2} = 0$, $rk =0$,
$sq =0$,
$dp =0$, $py =1$,  ${\rm vert}(\Delta^*) =9$:

\begin{tabular}{|c|c|c|c|c|c|c|c|c|}
-2& 0& 0& 0 & -1&-2& 1 &-2 &-1\\
 0& 1& 0& 0 &  1& 0& 0 &-1 &-1\\
 0& 0& 1& 0 & -1&-1& 0 & 0 & 1\\
-1& 0& 0& 1 &  0& 0& 0 & 0 & 0
\end{tabular}

$----------------------------------$

V(166):  $\deg = 36$,  $h^{1,2} = 0$, $rk =0$,
$sq =0$,
$dp =0$, $py =1$,  ${\rm vert}(\Delta^*) =9$:

\begin{tabular}{|c|c|c|c|c|c|c|c|c|}
-1& 0& 0& 0 & -2&-1&-1 &-2 & 1\\
 1& 1& 0& 0 & -1&-1& 1 & 0 & 0\\
 0& 0& 1& 0 &  0& 1&-1 &-1 & 0\\
-1& 0& 0& 1 &  0& 0& 0 & 0 & 0
\end{tabular}

$----------------------------------$

}

\end{document}